\documentclass[a4paper,10pt]{article}
\usepackage{anysize}
\marginsize{3.0cm}{3.0cm}{2.5cm}{1.5cm}

\usepackage{latexsym,amsfonts,amsmath}
\usepackage{graphicx}
\usepackage{epstopdf}

\newcommand{\widebar}[1]{\overline{#1}}

\newcommand{\defi}{\stackrel{\triangle}{=}}

\newtheorem{condition}{Condition}{\bfseries}{\itshape}

\newtheorem{theorem}{Theorem}{\bfseries}{\itshape}

\newtheorem{corollary}{Corollary}{\bfseries}{\itshape}

{\bfseries}{\itshape}

{\bfseries}{\itshape}

\newtheorem{lemma}{Lemma}{\bfseries}{\itshape}

\newtheorem{remark}{Remark}{\bfseries}{\itshape}

\newtheorem{definition}{Definition}{\bfseries}{\itshape}

\begin{document}
\author{
Alexey Piunovskiy\thanks{Corresponding author} \\
Department of Mathematical Sciences, University of Liverpool, L69 7ZL, UK.\\ \texttt{piunov@liv.ac.uk}
\and Yi Zhang \\
Department of Mathematical Sciences, University of Liverpool, L69 7ZL, UK. \\
\texttt{zy1985@liv.ac.uk} \\
}
\title{Discounted Continuous-time Markov Decision Processes with Unbounded Rates: the Dynamic Programming Approach}
\date{}

\maketitle

\par\noindent{\bf Abstract:} This paper deals with unconstrained discounted continuous-time Markov decision processes in Borel state and action spaces. Under some conditions imposed on the primitives, allowing unbounded transition rates and unbounded (from both above and below) cost rates, we show the regularity of the controlled process, which ensures the underlying models to be well defined. Then we develop the dynamic programming approach by showing that the Bellman equation is satisfied (by the optimal value). Finally, under some compactness-continuity conditions, we obtain the existence of a deterministic stationary optimal policy out of the class of randomized history-dependent policies. 

\section{Introduction}
In this paper, we show the existence of a deterministic stationary optimal policy out of
the class of randomized history-dependent policies for (unconstrained) discounted
continuous-time Markov decision processes (CTMDPs) with unbounded rates and with Borel state and action spaces. CTMDPs have been studied intensively since 1960s, and their formal constructions are available in \cite{Howard:1960} for deterministic stationary policies,
in \cite{Miller:1968} for deterministic Markov policies, and
in \cite{Hordijk:1979} for randomized Markov policies. The first
rigorous construction allowing deterministic
history-dependent policies is in \cite{Yushkevich:1977,Yushkevich:1980}, where the author viewed CTMDPs under deterministic history-dependent policies as special semi-Markov decision processes (SMDPs) whose actions are taken from spaces of measurable mappings. The first
successful construction of CTMDPs allowing randomized history-dependent policies is in \cite{Kitaev:1986}, which is based on
\cite{Jacod:1975}. As noted in \cite{Feinberg:2004}, although the
construction in \cite{Yushkevich:1977,Yushkevich:1980} is
restricted to deterministic history-dependent policies, it can be
modified to allow randomized history-dependent policies. In this
connection, Yushkevich's construction is indeed equivalent to
Kitaev's construction. To our best knowledge, currently, Kitaev's
construction provides the standard setup for CTMDPs allowing
randomized history-dependent policies, which we base the present work on. A brief reminder of this construction is provided below.

The expected total discounted cost has been a
common optimality criterion for CTMDPs optimization
problems\footnote{It is a standard practice to use ``CTMDPs'' and
``CTMDPs optimization problems'' interchangeably.}, and the
existence of an optimal policy for discounted CTMDPs has been
studied by numerous authors, see for example,
\cite{Feinberg:2004,Kakumanu:1971,Piunovskiy:1998,Yushkevich:1979}. In greater detail, \cite{Kakumanu:1971} is
restricted to deterministic Markov policies,
\cite{Yushkevich:1979} considers deterministic history-dependent
policies, while \cite{Feinberg:2004,Piunovskiy:1998} allow randomized history-dependent policies into
consideration. It should be emphasized that all of them assume
uniformly bounded transition rates. On the contrary,
\cite{GuoZhu:2002,Guo:2003} study discounted CTMDPs allowing transition rates to be not
uniformly bounded. However, the conditions
assumed therein are difficult for verifications, as some of them
are not directly imposed on the primitives but on the transition probability functions. Later on, there have
been developments in the direction of only imposing conditions on
the primitives, while still allowing unbounded transition rates,
see \cite{Guo:2007,Yan:2008} and the relevant chapters in the monograph \cite{Guo:2009}. It should be
noted that all of the aforementioned works allowing unbounded
transition rates are restricted to the class of randomized Markov
policies. As a fact of matter, according to \cite{Guo:2006TOP},
the study of CTMDPs with the combination of randomized
history-dependent policies and unbounded transition rates had been
an over thirty year-old open problem. To our best knowledge, the first
successful treatment for such CTMDPs is given by \cite{GuoABP:2010},
where the state space is countable.

In the present paper, we consider a more general case by allowing
randomized history-dependent policies, unbounded transition rates
and Borel state and action spaces into consideration, while all our conditions are imposed on
the primitives. The cost rates being allowed to be unbouned (both from below and above) are more general than those considered in \cite{GuoZhu:2002,Guo:2003,Guo:2003IEEE, Guo:2006TOP,Guo:2007,Guo:2009,GuoABP:2010} and many others, too.

The main contributions of the present paper are
triple-folded. Under the imposed conditions on the primitives, we firstly show the regularity of the controlled process under any given randomized history-dependent policy, which allows a formal optimization problem statement. Then we
develop the dynamic programming approach, by showing that the optimal value of the problem satisfies the corresponding Bellman equation. Finally, we establish the existence
of a deterministic stationary optimal policy. In relation to the most recent literature on this topic, the present work refines
\cite{Guo:2007} by considering randomized history-dependent policies\footnote{In comparison, \cite{Guo:2007} only considers a specific class of Markov policies, under which the resulting (nonhomogeneous) transition rates are required to be continuous in time, merely for the sake of validating the relevant results from \cite{Feller:1940}. In our opinion, this continuity is not needed.}, and extends
\cite{GuoABP:2010} to the case of Borel state spaces and more general cost rates.

The rest of this paper is organized as follows. In Section
\ref{MDPsectiondescription}, we briefly describe Kitaev's construction for CTMDPs, and present some preliminary results including the regularity, Kolmogorov's forward equations and Dynkin's formula for the controlled processes, which could be not Markov. In Section
\ref{MDPdyapproach}, we present the main statements. Section \ref{secexample} contains a new example. We finish this paper with a conclusion in Section \ref{MDPsecconclusion}. Several statements presented in this paper appeared without proofs in \cite{int15}.

\section{Preliminaries}\label{MDPsectiondescription}

The following denotations are frequently used throughout this
paper. $I$ stands for the indicator function. $\delta_{x}(\cdot)$
is the Dirac measure concentrated at $x.$ ${\cal{B}}(X)$ is the
Borel $\sigma$-algebra of the Borel space $X.$ ${\cal{F}}_1\bigvee{\cal{F}}_2$ is the smallest $\sigma$-algebra
containing the two $\sigma$-algebras ${\cal{F}}_1$ and
${\cal{F}}_2.$ $\mathbb{R}_+\defi (0,\infty).$
$\mathbb{R}_+^0\defi[0,\infty)$. $\mathbb{Z}_+^0\defi\mathbb{N}\bigcup\{0\}$. The abbreviation $s.t.$ (resp.
$a.s.$) stands for ``subject to'' (resp. ``almost surely'').

\subsection{Kitaev's construction}
The materials presented in this subsection are mainly from
\cite{Kitaev:1986,Kitaev:1995,Piunovskiy:1998}.

The primitives
of discounted CTMDPs are the following elements:
\begin{itemize}
\item state space: $(S,{\cal{B}}(S))$ (arbitrary Borel), \item action space:
$(A,{\cal{B}}(A))$ (arbitrary Borel),
\item admissible action space $A(x)\in {\cal
B}(A)$ and the space of admissible action-state pairs
$K\defi\{(x,a)\in S\times A: a\in A(x)\}\in {\cal{B}}(S\times A),$
assumed to contain the graph of a measurable function $\phi$ from
$S$ to $A$ such that $\forall~x\in S,~\phi(x)\in A(x),$
\item
transition rate: $q(dy|x,a),$ a signed kernel on ${\cal{B}}(S)$
given $(x,a)\in K$, taking nonnegative values on
$\Gamma_S\setminus\{x\}$ with $\Gamma_S\in{\cal{B}}(S),$ being
conservative in the sense of $q(S|x,a)=0$ and stable in that
$\bar{q}_x=\sup_{a\in A(x)}q_x(a)<\infty,$ where $q_x(a)\defi
-q(\{x\}|x,a),$
\item cost rate: $c_0(x,a)$ measurable in
$(x,a)\in K,$
\item discount factor:
$\alpha>0,$
\item initial distribution: $\gamma(\cdot),$ a
probability measure on $(S,{\cal B}(S)).$
\end{itemize}
Incidentally, we remind that a singleton $\{x\}\subseteq
S$ is measurable, and $q_x(a)$ is measurable on $K$, see
\cite[Prop 7.29]{Bertsekas:1978}. In what follows, for the sake of
formality, if needed, $\forall~\Gamma_S\in{\cal B}(S),$ we may
consider $q(\Gamma_S|x,a)$ as its measurable extension on $S\times
A,$ where $q(\Gamma_S|x,a)=0$ on $(S\times A)\setminus K,$ and
similar assertions are applicable to other functions such as
$c_0,$ and so on. This is just the convention, see
\cite[Chap.6]{Hernandez-Lerma:1996}.

Given the above primitives, let us recall the
construction of the underlying stochastic basis $(\Omega,{\cal
F},\{{\cal F}_t\}_{t\ge 0},P_\gamma^\pi)$ and the controlled
process $\{\xi_t,t\ge 0\}$ thereon, as given in \cite{Kitaev:1986} (see also \cite{Kitaev:1995,Piunovskiy:1998} for more details). This is done in four steps.

Step 1: measurable space $(\Omega,{\cal{F}}).$ Having firstly defined the measurable space of $(\Omega^0,{\cal{F}}^0)\defi ((S\times
\mathbb{R}_+)^\infty,{\cal{B}}((S\times \mathbb{R}_+)^\infty)),$
let us adjoin all the sequences of the form
$$(x_0,~\theta_1,~x_1,~\dots,~\theta_{m-1},~x_{m-1},~\infty,~x_{\infty}~,\infty,~x_{\infty},~\dots)$$
to $\Omega^0,$ where $x_l\in S$, $x_\infty\notin S$ is
an isolated point, $ m\ge 1$ is some integer, $\theta_l\in\mathbb{R}_+$ and $x_l\ne x_\infty$
for all nonnegative integers $l\le m-1$. After the corresponding
modification of the $\sigma$-algebra ${\cal{F}}^0,$ we obtain the
basic measurable space $(\Omega, {\cal{F}}).$

Step 2: stochastic process $\{\xi_t,t\ge 0\}$ and history $\{{\cal
F}_t\}_{t\ge 0}.$ Putting $T_0\defi
0,~T_m\defi\theta_1+\theta_2+\dots+\theta_m$,
$T_\infty\defi\lim_{m\rightarrow \infty}T_m,$ we can define the
process of interest:
$$\xi_t(\omega)\defi\sum_{m\ge 0}I\{T_m\le
t<T_{m+1}\}x_m+I\{T_\infty\le t\}x_\infty$$ together with the
history it is adapted to:
$${\cal{F}}_t\defi \sigma(\{T_m\le s, x_m\in \Gamma_S\}:
\Gamma_S\in {\cal{B}}(S),s\le t,m\ge 0).$$
In what follows, as usual, $\omega=\{x_0,\theta_1,x_1,\dots\}$ is often
omitted, and $h_m(\omega)=(x_0,\theta_1,\dots, \theta_m,x_m)$ is referred
to as an $m$-component history. Here, $\theta_m$ (resp. $T_m$,
$x_m$) can be understood as the sojourn
times (resp. the jump moments, the state of the process on the
interval $[T_m,T_{m+1})$). We do not intend to consider the
process after $T_\infty:$ the isolated point $x_\infty$ will be
regarded as absorbing.

Step 3: policy $\pi.$  Having adjoint the isolated point
$a_\infty$ to $A$, we thus define $A_\infty\defi A\bigcup
\{a_\infty\},$ and put $A(x_\infty)\defi \{a_\infty\}$. Similarly,
$S_\infty\defi S\bigcup \{s_\infty\}$. Denoting ${\cal{F}}_{s-}\defi
\bigvee_{t<s}{\cal{F}}_t,$ the predictable (with respect to
$\{{\cal{F}}_t\}_{t\ge0}$) $\sigma$-algebra ${\cal{P}}$ on
$\Omega\times \mathbb{R}_+^0$ is given by
${\cal{P}}\defi \sigma(\Gamma\times \{0\} ~(\Gamma\in {\cal{F}}_0), \Gamma\times
(s,\infty)~(\Gamma\in {\cal{F}}_{s-})).$ See
\cite[Chap.4]{Kitaev:1995} for more details. Now the following definitions are in position:
\begin{itemize}
 \item Randomized history-dependent policy: $\pi(\cdot|\omega,t)$,
a ${\cal{P}}$-measurable transition probability function on
$(A_\infty, {\cal{B}}(A_\infty))$, concentrated on $A(\xi_{t-}).$ Below, $U$ is the set of all such policies.

\item Randomized Markov policy: $\pi(\cdot|\omega,t)=\pi^m(\cdot|\xi_{t-}(\omega),t).$
Here concerning the RHS, $\pi^m(\cdot|x,t)$ is ${\cal
B}(S_\infty\times \mathbb{R}_+^0)$-measurable.

\item Randomized stationary policy:
$\pi(\cdot|\omega,t)=\pi^s(\cdot|\xi_{t-}(\omega)).$ Here concerning
the RHS, $\pi^s(\cdot|x)$ is ${\cal B}(S_\infty)$-measurable.

\item Deterministic stationary policy: $\pi(\cdot|\omega,t)=I\{\cdot\ni \phi(\xi_{t-}(\omega)) \},$ where $\phi: S_\infty\rightarrow A_\infty$ is a measurable mapping. Such policies are denoted as $\phi.$
\end{itemize}

\begin{remark}
The term ``randomized policies'' is adopted from \cite{Feinberg:2004,Kitaev:1986,Piunovskiy:1998}. However, under a randomized policy, it does not mean that decisions are made randomly continuously in time, which is not always possible (see \cite[Sec.7]{Feinberg:2004}). In fact, the term of randomized policies should be understood as relaxed control policies, as remarked in \cite[Chap.4]{Kitaev:1995}. Throughout this paper, the most general policy under consideration is randomized history-dependent.
\end{remark}

Step 4: ($\gamma,\pi$-dependent) probability measure
$P_{\gamma}^\pi$ on $(\Omega,{\cal{F}}).$ Under any fixed policy
$\pi$, let us define
\begin{eqnarray}\label{MDPrandommeasure2}
\nu^\pi(\omega,\Gamma_S\times dt)\defi
\Lambda(\Gamma_S|\omega,t)dt\defi\left[\int_{A}\pi(da|\omega,t)q(\Gamma_S\setminus\{\xi_{t_-}\}|\xi_{t-},a)\right]dt,
\end{eqnarray} where $\Gamma_S\in{\cal{B}}(S),$ and the obvious dependence of
$\Lambda$ on $\pi$ has been omitted. This random measure is
predictable, see \cite{Kitaev:1986,Kitaev:1995,Piunovskiy:1998}.
According to \cite[Chap.4]{Kitaev:1995} (see also
\cite{Jacod:1975}), the ``jump intensity'' $\Lambda$ has the following form:
\begin{eqnarray}\label{MDPconstructintensity}
\Lambda(dy|\omega,t)&=&\sum_{m\ge 0}I\{T_m<t\le
T_{m+1}\}\Lambda^m(dy|x_0,\theta_1,\dots,x_m,t-T_m)\nonumber\\
&&+I\{t=0\}\Lambda^0(dy|x_0),
\end{eqnarray}
where $\forall~\Gamma_S\in {\cal{B}}(S)$,
$\Lambda^m(\Gamma_S|x_0,\theta_1,\dots,x_m,u)$ are some
nonnegative, non-random measurable functions. Then comparing
(\ref{MDPrandommeasure2}) with (\ref{MDPconstructintensity}), we
have the explicit formula\footnote{In fact, since $\pi(\cdot|\omega,t)$ is ${\cal P}$-measurable, it also admits a similar representation to $\Lambda(\cdot|\omega,t)$ (see (\ref{MDPconstructintensity})). This is because of \cite[Chap.4]{Kitaev:1995}. In this connectation, to be absolutely rigorous, one should write $\pi^m(\cdot|x_0,\theta_1,\dots,x_m,u)$ in (\ref{MDPexplicitintensity}), rather than $\pi(\cdot|x_0,\theta_1,\dots,x_m,u+T_m).$ Nevertheless, here and below, we omit that superscript $m$, and use the denotation $\pi(\cdot|x_0,\theta_1,\dots,x_m,u+T_m)$ for $\pi^m(\cdot|x_0,\theta_1,\dots,x_m,u).$ This is merely for brevity, as the context always excludes any confusion; besides, the superscript $m$ has already been used to indicate a Markov policy.} for $\Lambda^m:$
\begin{eqnarray}\label{MDPexplicitintensity}
\Lambda^m(dy|x_0,\theta_1,\dots,x_m,u)
=\int_{A}\pi(da|x_0,\theta_1,\dots,x_m,u+T_m)q(dy\setminus\{x_m\}|x_m,a).
\end{eqnarray}
Let $\hat{H}_0\defi S$ and $\hat{H}_m\defi S\times
((0,\infty]\times S_\infty)^m,m=1,2,\dots.$ The marginal of $P_\gamma^\pi$ on
$\hat{H}_0$ coincides with $\gamma.$\footnote{Below, with some abuse of denotation,
we also use $P_\gamma^\pi$ for the marginals on $\hat{H}_m$.} Suppose that
$P_\gamma^\pi$ on $\hat{H}_m$ for $1\le m\le k$ has been constructed. Now it is only needed to
construct $P_\gamma^\pi$ on $\hat{H}_{k+1}.$ But this can be done
via
\begin{eqnarray}\label{MDPconstructprob}
P^\pi_\gamma(\Gamma^{\hat{H}_k}\times (du\times dy))&\defi&
\int_{\Gamma^{\hat{H}_k}}P_\gamma^{\pi}(dh_k)I\{\theta_k<\infty\}\Lambda^k(dy|h_k,u) e^{-\int_0^u\Lambda^k(S|h_k,v)dv}du, \nonumber\\
&&\nonumber \\
P^{\pi}_\gamma(\Gamma^{\hat{H}_k}\times
(\infty,x_\infty))&\defi&
\int_{\Gamma^{\hat{H}_k}}P_{\gamma}^{\pi}(dh_k)\left\{I\{\theta_k=\infty\}+I\{\theta_k<\infty\}e^{-\int_{0}^\infty\Lambda^k(S|h_k,v)dv}\right\},
\end{eqnarray}
where $\Gamma^{\hat{H}_k}\in{{\cal{B}}(\hat{H}_k)}.$ It
remains to apply the induction and Ionescu-Tulcea's theorem
\cite[p.140-141, Prop.7.28]{Bertsekas:1978} to induce that
$P^\pi_\gamma$ is the unique probability measure on
$(\Omega,{\cal{F}})$ such that its projection (marginal) onto
$\hat{H}_m$ satisfies (\ref{MDPconstructprob}), $m=0,1,\dots.$
This gives rise to stochastic basis $(\Omega, {\cal F},\{{\cal F}_t\}_{t\ge
0}, P_\gamma^\pi),$ which is always assumed to be complete.

In fact, according to \cite {Kitaev:1986}, if we define the random
measure
\begin{equation}\label{MDPrandommeasure1}
\mu(\omega,dt,dy)\defi\sum_{m\ge 1}I\{T_m<\infty\}I\{x_m\in
dy\}I\{T_m\in dt\},
\end{equation}
then under any fixed policy $\pi$ and initial
distribution $\gamma$, the above defined measure $P^\pi_{\gamma}$ on
$(\Omega, {\cal{F}})$ is such that its projection on the
$0$-component history is $\gamma,$ and $\nu^\pi$ defined by
(\ref{MDPrandommeasure2}) is the dual predictable projection of
$\mu$ defined by (\ref{MDPrandommeasure1}). See \cite[Chap.4]{Kitaev:1995}
for more details.

Below, when $\gamma(\cdot)$ is a Dirac measure
concentrated at $x\in S,$ we use the ``degenerated''
denotation $P_x^\pi.$ Expectations with respect to
$P_\gamma^\pi$ and $P_x^\pi$ are denoted as $E_{\gamma}^\pi$ and
$E_{x}^\pi,$ respectively.

\subsection{Properties of the controlled process and optimization problem statement}
\begin{condition}\label{MDPRegularitycons}
There exist a measurable (weight) function $w(x)\ge 1$ on $S$ and
constants $\rho\ge0,b\ge 0$ such that
\par\noindent (a) $\bigcup_{l=0}^{\infty}S_{l}=S$ and $\lim_{l\rightarrow\infty}\inf_{x\in S\setminus
S_{l}}w(x)=\infty$ for an increasing system of measurable subsets $S_{l}\subseteq S.$
\par\noindent (b) $\int_{S}q(dy|x,a)w(y)\le \rho w(x)+b,\forall~x\in S, a\in A(x).$
\par\noindent (c) For any $l\in \mathbb{Z}_+^0,$ $\sup_{x\in
{S_l}}\bar{q}_x<\infty,$ where $S_l$ has been defined in part (a),
and $\bar{q}_{x}\defi\sup_{a\in A(x)}q_{x}(a).$
\end{condition}

\begin{remark}\label{MDPrho} Below, we assume $\rho>0,$ where
$\rho$ is defined in Condition \ref{MDPRegularitycons}. This can be done without loss of generality, because
the case of $\rho=0$ can always be considered by passing to the
limit as $\hat{\rho}\rightarrow 0,$ with $\hat{\rho}>0.$ We emphasize that if Condition \ref{MDPRegularitycons} is satsified by $\rho=0,$ it is also satisfied by any arbitrarily fixed $\hat{\rho}>0.$
\end{remark}
Condition \ref{MDPRegularitycons} is of a Lyapunov type. Theorem \ref{MDPregularitythm} shows that it guarantees the $\xi_t$ process to be non-explosive.
\begin{condition}\label{MDPfinitenesscons}
\par\noindent (a) $\int_{S}\gamma(dy)w(y)<\infty,$ where
$\gamma$ is the given initial distribution.
\par\noindent (b) $\alpha>\rho,$ where $\alpha$ is the discount
factor, and $\rho$ is as in Condition \ref{MDPRegularitycons}.
\par\noindent (c) There exist constants $M\ge 0$ and $c\ge 0$ such that
$|\inf_{a\in A(x)}c_0(x,a)|\le Mw(x)+c,\forall~x\in S.$
\end{condition}
This condition guarantees that the performance functional (\ref{MDPperformancefunctional2}) is well defined. Condition \ref{MDPfinitenesscons}(c) is a version of the one imposed in \cite{Piunovskiy:1989}, where the author studies CTMDPs with bounded transition rates and average criteria.

\begin{theorem}\label{MDPregularitythm}
Suppose Condition \ref{MDPRegularitycons} is satisfied. Then under any policy
$\pi\in U$, the following assertions hold:
\par\noindent (a)  For any given initial distribution $\gamma,$
$P_{\gamma}^\pi(T_\infty=\infty)=1,$ and hence $\forall~t\ge 0$, $P^{\pi}_\gamma(\xi_t\in
S)=1.$ So explosion does not occur. Moreover, for all $x\in S,t\ge 0$,
\begin{eqnarray*}
E_{x}^\pi\left[w(\xi_{t})\right]\le e^{\rho t}w(x)+\frac{b}{\rho}(e^{\rho t}-1).
\end{eqnarray*}

\par\noindent (b) If additionally Condition
\ref{MDPfinitenesscons} is satisfied, then for any $\gamma,$ inequality
$$
V_0(\pi)\ge-\frac{M(\alpha
\int_S\gamma(dy)w(y)+b)}{\alpha(\alpha-\rho)}-\frac{c}{\alpha}>-\infty
$$
holds, where
\begin{eqnarray}\label{MDPperformancefunctional2}
V_0(\pi)\defi E_{\gamma}^\pi\left[\int_0^\infty e^{-\alpha
t}\int_{A}c_0(\xi_{t-},a)\pi(da|\omega,t)dt\right].
\end{eqnarray}
\end{theorem}

We use denotation $V_0(x,\pi)$ if the initial distribution $\gamma$ is concentrated at state $x\in S$.

The proofs of this theorem and the other main statements presented
in this paper can be found in the appendix.

Theorem \ref{MDPregularitythm} implies that the following CTMDPs optimization problem under consideration is well defined:
\begin{eqnarray}\label{MDPprobstatement}
&&V_0(\pi)\rightarrow \min_{\pi\in U}.
\end{eqnarray}

\begin{definition}
Denote by $V_0^\ast\defi \inf_{\pi\in
U}V_0(\pi)$ the optimal value of CTMDP (\ref{MDPprobstatement}). A policy $\pi^\ast$ is called optimal, if
$V_0(\pi^\ast)=V_0^\ast.$ CTMDP (\ref{MDPprobstatement}) is called
solvable, if such a $\pi^\ast$ exists.
\end{definition}
\begin{remark}\label{MDPfeasibleassump} Equality (\ref{MDPexplicitintensity}) holds $P^\pi_\gamma$-a.s.,
as well as all the subsequent equalities and inequalities
involving $\omega$. The values of integrals like (\ref{MDPperformancefunctional2}) do not change, if we replace $\xi_{t-}$ with $\xi_t.$
\end{remark}

\subsection{Auxiliary results}
Generally speaking, $\bar{q}_x$ may be not measurable. However, according to \cite[D.5 Prop.]{Hernandez-Lerma:1996} (see also \cite[Prop.7.33]{Bertsekas:1978}), $\bar{q}_x$ is measurable on $S$ if  the following condition is satisfied. 
\begin{condition}\label{MDPmeas1}
\par\noindent(a) $A(x)$ is compact, $\forall~x\in S$.
\par\noindent(b) $q_x(a)$ is upper semicontinuous on $A(x),$ $\forall~x\in S.$
\end{condition} 

Kolmogorov's forward equation (in the integral form) and Dynkin's formula are rather useful tools for studying CTMDPs. In case $\pi$ is Markov, they are well known. For a randomized history-dependent policy $\pi,$ under the imposed conditions, it turns out that they still hold.

\begin{condition}\label{MDPstrongregularitycons}
There exists a constant $L>0$ such that $0\le \bar{q}_{x}<Lw(x)$, $\forall~x\in S.$
\end{condition}
We need this condition to be sure that the last term in formula (\ref{MDPforwardkolmogoroveqn3}) is finite.

\begin{theorem}
\label{MDPkolmogorov}
\par\noindent(a)
Suppose Condition \ref{MDPRegularitycons} is
satisfied. Then under any fixed policy $\pi$, $\forall~x\in S,t\in \mathbb{R}_+^0$,
$\forall~\Gamma\in {\cal{B}}(S)$ such
that $\exists~l:$ $\Gamma\subseteq S_l,$ with $S_l$ being defined in
Condition \ref{MDPRegularitycons}, Kolmogorov's forward equation (in
the integral form) holds:
\begin{eqnarray}\label{MDPforwardkolmogoroveqn2}
P_{x}^\pi(\xi_t\in
\Gamma)&=&I\{x\in\Gamma\}+E_{x}^\pi\left[\int_0^t\int_A\pi(da|\omega,u)q(\Gamma\setminus\{\xi_{u}\}|\xi_{u},a)du\right]\nonumber\\
&&-E_{x}^\pi\left[\int_{0}^t\int_A\pi(da|\omega,u)q_{\xi_{u}}(a)I\{\xi_{u}\in\Gamma\}
du\right].
\end{eqnarray}

\par\noindent(b) In part (a), if we replace Condition \ref{MDPRegularitycons}(c) by
Condition \ref{MDPstrongregularitycons}, whereas all the other parts of Condition \ref{MDPRegularitycons} are still satisfied, then we have the following stronger statement: $\forall~\Gamma \in
{\cal{B}}(S),$
\begin{eqnarray}\label{MDPforwardkolmogoroveqn3}
P_{x}^\pi(\xi_t\in
\Gamma)&=&I\{x\in\Gamma\}+E_{x}^\pi\left[\int_0^t\int_A\pi(da|\omega,u)q(\Gamma\setminus\{\xi_{u}\}|\xi_{u},a)du\right]\nonumber\\
&&-E_{x}^\pi\left[\int_{0}^t\int_A\pi(da|\omega,u)q_{\xi_{u}}(a)I\{\xi_{u}\in\Gamma\}
du\right].
\end{eqnarray}

The expectations that appear in the above formulae are finite.
\end{theorem}

For the case of uniformly bounded
$\bar{q}_x$, Kolmogorov's forward equation (\ref{MDPforwardkolmogoroveqn3}) has been established
in \cite[Lem.4]{Kitaev:1986}. Throughout this paper, Condition \ref{MDPstrongregularitycons} is only required for proving Theorem \ref{MDPkolmogorov}(b), while Theorem \ref{MDPkolmogorov}(b) itself is never used elsewhere in this paper. However, it is needed in \cite{ABPZY:20102}.

We need parts (a,b) of the next condition for establishing Dynkin's formula, where the product $\widebar q_{\xi_v} u(\xi_v)$ must be integrable for $u\in {\textbf B}_{w'}(S)$. (See Definition \ref{defin2}.)
\begin{condition}\label{MDPdpcons}
There exist a measurable function $w'(x)\ge 1$ on $S$ and
nonnegative constants $L',\rho'$ and $b'$ such that the following
assertions hold:
\par\noindent (a) $(\bar{q}_x+1)w'(x)\le L'w(x),$
where $w$ comes from Condition \ref{MDPRegularitycons}.
\par\noindent (b) $\int_{S}q(dy|x,a)w'(y)\le \rho' w'(x)+b',\forall~x\in S, a\in A(x).$
\par\noindent (c) $\alpha>\rho'.$
\par\noindent (d) There exist constants $M'\ge0$ and $c'\ge 0$ satisfying $|\inf_{a\in A(x)}c_0(x,a)|\le M'w'(x)+c',\forall~x\in S.$
\end{condition}

Condition \ref{MDPdpcons}(c,d) guarantees that the corresponding performance functional is well defined (cf Condition \ref{MDPfinitenesscons}(b,c) ). Under Condition \ref{MDPRegularitycons} and Condition \ref{MDPdpcons}(a), $E^\pi_x[w'(\xi_t)]<\infty$ due to Theorem \ref{MDPregularitythm}(a).

\begin{definition}\label{defin2} A measurable function $u$ on $S$ satisfying
$\sup_{x\in S}\frac{|u(x)|}{w(x)}<\infty$ (resp. $\sup_{x\in
S}\frac{|u(x)|}{w'(x)}<\infty$) is said to have a bounded $w$-(resp. $w'$-)weighted norm, with the norm $||u||_w\defi
\sup_{x\in S}\frac{|u(x)|}{w(x)}$ (resp. $||u||_{w'}\defi
\sup_{x\in S}\frac{|u(x)|}{w'(x)}$). The collection of all
functions $u$ on $S$ with a bounded $w$-(resp. $w'$-)weighted
norm is denoted by $\textbf{B}_w(S)$ (resp.
$\textbf{B}_{w'}(S)$).
\end{definition}

\begin{theorem}
\label{MDPDynkinsformula} Suppose Condition
\ref{MDPRegularitycons} and Condition
\ref{MDPdpcons}(a,b) are satisfied. Then $\forall~u\in \textbf{B}_{w'}(S),$ the following two versions of Dynkin's
formula hold:
\begin{eqnarray}\label{MDPDynkinformula1}
E_x^{\pi}[u(\xi_t)]-u(x)=E_x^\pi\left[\int_0^t\int_S\int_A\pi(da|\omega,v)q(dy|\xi_{v},a)u(y)dv\right],
\end{eqnarray}
\begin{eqnarray}\label{MDPDynkinformula2}
E_x^\pi[u(\xi_t)]e^{-\alpha t}-u(x)=E_x^\pi\left[\int_0^t
e^{-\alpha v}\left\{-\alpha
u(\xi_v)+\int_S\int_A\pi(da|\omega,v)q(dy|\xi_v,a)u(y)\right\}dv\right].
\end{eqnarray}
\end{theorem}

\section{Main statements}\label{MDPdyapproach}
\begin{condition}\label{MDPmeas2}
\par\noindent (a) For any bounded nonnegative
measurable function $u(y)$ on $S$ and fixed $x\in S$, $u'(x,a)\defi\int_Su(y)q(dy|x,a)$ is
lower semicontinuous in $a\in A(x).$
\par\noindent (b) $\int_S w(y)q(dy|x,a)$ is continuous in $a\in
A(x)$, $\forall~x\in S,$ where $w$ comes from Condition
\ref{MDPRegularitycons}.
\par\noindent (c) $c_0(x,a)$ is lower semicontinuous in $a\in A(x),
\forall~x\in S.$
\par\noindent (d) $A(x)$ is compact, $\forall~x\in S$.
\end{condition}
\begin{remark}\label{MDPremarkextra}
By reasoning similarly to \cite[p.44]{Hernandez-Lerma:1999},  one can show that Condition \ref{MDPmeas2}(a) is equivalent to the following: for any $x\in S$ and bounded measurable function $u(y)$
on $S,$ function $\int_Su(y)q(dy|x,a)$ is continuous in $a\in A(x)$. Therefore, Condition \ref{MDPmeas2}(a) is stronger than
Condition \ref{MDPmeas1}(b).
\end{remark}

The next statement is similar to Theorem 3.3 (b) in \cite{Guo:2007}.

\begin{theorem}\label{MDPguosprop}
Suppose Condition \ref{MDPRegularitycons}(b), Condition \ref{MDPfinitenesscons}(b,c) and Condition
\ref{MDPmeas2} are satisfied. Then the Bellman equation
\begin{equation}\label{MDPdpeqn}
\alpha u(x)=\inf_{a\in
A(x)}\left\{c_0(x,a)+\int_Sq(dy|x,a)u(y)\right\}.
\end{equation}
admits a
solution $u^\ast\in \textbf{B}_w(S),$ which is given by the point-wise limit of
the following non-increasing sequence of measurable functions
$\{u^{(n)},n=0,1,\dots\}$:
\begin{eqnarray}\label{eqm}
u^{(0)}(x)&\defi& \frac{M(\alpha
w(x)+b)}{\alpha(\alpha-\rho)}+\frac{c}{\alpha},\nonumber \\
u^{(n+1)}(x)&\defi&
\inf_{a\in
A(x)}\left\{\frac{c_0(x,a)}{\alpha+1+\bar{q}_x}+\frac{1+\bar{q}_x}{\alpha+1+\bar{q}_x}\int_S
u^{(n)}(y)\left(\frac{q(dy|x,a)}{1+\bar{q}_x}+I\{x\in dy\}\right)\right\}.
\end{eqnarray}
For each $n=0,1,2,\ldots$
  $$|u^{(n)}(x)|\le \frac{M(\alpha
w(x)+b)}{\alpha(\alpha-\rho)}+\frac{c}{\alpha}.$$
\end{theorem}

\begin{remark}\label{MDPdyremark} (a) Suppose Condition \ref{MDPdpcons}(b,c,d) is satisfied. If additionally Condition
\ref{MDPmeas2} (with $w$ being replaced with $w'$ in its part (b)) is satisfied, then the statements of Theorem \ref{MDPguosprop} are still valid, with
$w,M,c,\rho$ and $b$ being replaced by $w',M',c',\rho'$ and $b'$
everywhere. This remark can be verified by repeating the reasonings used in the proof of Theorem \ref{MDPguosprop}, with obvious modifications.

(b) Condition \ref{MDPfinitenesscons}(b), Condition \ref{MDPdpcons}(a) and Condition \ref{MDPmeas2} altogether imply that $\int_S w'(y) q(dy|x,a)$ is continuous in $a\in A(x)$ for each $x\in S$ (see \cite[Lem.8.3.7.]{Hernandez-Lerma:1999}).
\end{remark}

\begin{theorem}\label{MDPbellmanthm}
Suppose Condition \ref{MDPRegularitycons},
Condition \ref{MDPfinitenesscons}(a,b), Condition \ref{MDPdpcons} and
Condition \ref{MDPmeas2} are
satisfied. Then the following assertions hold:
\par\noindent (a) Suppose function $u^\ast\in\textbf{B}_{w'}(S)$ solves the Bellman equation
(\ref{MDPdpeqn}), then, for some deterministic stationary policy $\phi^*$
$$\int_S\gamma(dy)u^\ast(y)=\inf_{\pi}V_0(\pi)=V_0(\phi^*).$$
If a measurable map $\phi^*:~x\to\phi^*(x)\in A(x)$ provides the infimum in (\ref{MDPdpeqn}) then policy $\phi^*$ is optimal.
\par\noindent (b) The Bellman equation (\ref{MDPdpeqn}) has a unique solution $u^*$ in the class $\textbf{B}_{w'}(S)$ which can be constructed using iterations (\ref{eqm}), where $w,M,c,\rho$ and $b$ should be replaced with $w',M',c',\rho'$ and $b'$.
\par\noindent (c) The Bellman function $u^\ast$ solves the following dual linear program (DLP) in the space of measurable functions on $S$:
\begin{eqnarray}\label{MDPunconstrainedDLP}
&&\int_{S}\gamma(dy)v(y)\rightarrow \max_v\\
&s.t.&\nonumber\\
&&\frac{1}{\alpha}c_0(x,a)-v(x)+\frac{1}{\alpha}\int_Sv(y)q(dy|x,a)\ge
0, \forall~(x,a)\in K;\nonumber\\
&&v\in \textbf{B}_{w'}(S).\nonumber
\end{eqnarray}
\par\noindent (d) Suppose $v$ is feasible for DLP (\ref{MDPunconstrainedDLP}). Then it solves the DLP if and only if $v(x)=u^\ast(x)$ a.s. (with respect to $\gamma$).
\end{theorem}

\section{Example} \label{secexample}
Consider a one-channel queuing system without any space for waiting: any job that finds the server busy is rejected. We characterize every job by its volume $x\in(0,1]$, so that the state space is $S=[0,1]$: $\xi_t=0$ means the system is idle; $\xi_t=x\in(0,1]$ means the corresponding job is under service. We put $A=[0,\infty)$, and action $a\in A$ represents the service intensity. Let $A(0)=0$ and $A(x)=\left[0,\frac{\bar A}{x}\right],$ where $\bar{A}\ge0$ is a constant. The jobs arrive according to a Poisson process with a fixed rate $\lambda>0$, and the volume is distributed according to density $5x^4$, $x\in(0,1]$ independently of anything else. Therefore,
  $$q(\Gamma|0,a)=5\lambda\int_{\Gamma\setminus\{0\}}y^4 dy-\lambda I\{\Gamma\ni 0\},\forall~\Gamma\in {\cal B}([0,1]).$$
For any fixed $x\in(0,1],a\in A(x)$, the service time of a job of volume $x$ is exponentially distributed with parameter $\frac{a}{x}$, so that
  $$q(\Gamma|x,a)=I\{0\in\Gamma, x\notin\Gamma\}\frac{a}{x}-I\{0\notin\Gamma, x\in\Gamma\}\frac{a}{x},\forall~\Gamma\in{\cal B}([0,1]).$$
  We assume that when a served job leaves the system, it gives an income of one unit; the holding cost of a job of volume $x\in(0,1]$ equals $C_1x$ per time unit; and the service intensity $a\in A$ is associated with the cost rate $C_2a^2$. Here $C_1\ge0$ and $C_2\ge0$ are two constants. Thus
    $$c_0(x,a)=C_1x+C_2a^2-\frac{a}{x},\forall~x\in(0,1],a\in A(x),$$
and $c_0(0,0)=0$. We emphasize that as can be easily verified, $\bar{q}_x$ is unbouned, and $c_0(x,a)$ is unbouned (from both above and below) when $\bar{A}>\frac{1}{C_2}.$

Finally, let $\alpha$, the discount factor, be big enough:
  $$\alpha>4\lambda,$$
and let $\gamma$, the initial distribution, be such that
  $$\int_0^1\gamma(dy)\frac{1}{y^4}<\infty.$$

\begin{theorem}\label{t6} (a) For the model described, all the conditions formulated in this paper are satisfied. \\
\par\noindent(b) Suppose $C_1\ge 0$ is small enough (or $\alpha$ is big) in that $\frac{C_1}{2\alpha}\le 1,$ and define
\begin{eqnarray}\label{MDPexamplefunction1}
u(x,z)\defi-2\alpha C_2 x^2-z+2\sqrt{\alpha^2 C_2^2 x^4+C_1C_2 x^3+\alpha C_2 x^2 z}, \forall~x\in(0,1],z\in[0,\infty).
\end{eqnarray}
Then the following recursion relations
\begin{eqnarray*}
z^{(0)} & = & 0; \nonumber\\
u^{(n)}(x) & = & u(x,z^{(n)})= -2\alpha C_2 x^2-z^{(n)} +2\sqrt{\alpha^2 C_2^2 x^4+C_1C_2 x^3+\alpha C_2 x^2 z^{(n)}},~~~ x\in(0,1];\nonumber \\
z^{(n)} & = & 1-\frac{5\lambda}{\alpha+\lambda}\int_0^1 u^{(n)}(y) y^4 dy,~~~~~n=0,1,2,\ldots \nonumber
\end{eqnarray*}
converge: the sequence $\{z^{(n)},n=0,1,\dots\}$ is increasing and has a finite limit $z^*=\lim_{n\to\infty} z^{(n)}$, and $\lim_{n\rightarrow\ \infty}u^{(n)}(x)=u(x,z^\ast)\defi u^*(x), \forall~x\in (0,1].$ \\
\par\noindent (c) Suppose $\frac{C_1}{2\alpha}\le 1,$ and constant $\bar A$ is big enough in that the limiting function $u^*(x)$ satisfies inequality $\frac{u^*(x)+z^*}{2C_2}\le \bar A,\forall~x\in(0,1].$ Then $u^*(x)$, supplemented at zero by the value $u^*(0)\defi 1-z^*$, solves the Bellman equation (\ref{MDPdpeqn}), and the deterministic stationary policy
\begin{eqnarray}\label{MDPexamplepolicy}
\phi^*(x)=\frac{u^*(x)+z^*}{2xC_2},\forall~x\in(0,1],\mbox{ and } \phi^\ast(0)=0
\end{eqnarray}
is optimal.
\end{theorem}

\begin{remark}\label{rem6}
\par\noindent(a) If parameter $\bar A$ increases, the solution to this example does not change. We cannot put $A(x)=[0,\infty)$ because in this case the transition rate becomes unstable: $\sup_{a\in A(x)} q_x(a)=+\infty$.
\par\noindent(b) It follows from the proof of Theorem \ref{t6} that $z^*<\frac{10}{7} C_2\alpha+\frac{\alpha+\lambda}{\alpha}$ and function $u(x,z)$ defined by (\ref{MDPexamplefunction1}) decreases with $z$ for any fixed $x\in (0,1]$. These observations allow us to estimate the admissible values of $\bar A$.
\par\noindent(c) In case $C_1$ is very big (see part (c) of Theorem \ref{t6}) then it can happen that action $a^*=0$ becomes optimal for small values of $\xi_t=x$. Indeed, if $a>0$ then there can be transitions $x\to 0\to y\to\ldots$ with a good chance to have a big value of $y$ leading to a big holding cost in the future. Thus, in this situation it can be reasonable to select $a^*=0$ and finish with the cost rate $\frac{C_1 x}{\alpha},$ which is small if $x$ is small.
\end{remark}

\section{Conclusion}\label{MDPsecconclusion}
As mentioned in \cite{Hu:2003}, the standard results for (unconstrained) discounted CTMDPs include that the model is  well defined, the Bellman equation is satisfied, and there exists a deterministic stationary optimal policy. In the present work, taking into account as general as randomized history-dependent policies, we obtain all such standard results for CTMDPs in Borel spaces. The conditions we base our study on are imposed on the primitives, allowing unbounded transition and cost rates. In particular, our conditions imposed on the cost rate are more general than those in all the papers on discounted CTMDPs in the references. In this connection, the present paper is arguably in quite a general setup.

We emphasize that our conditions are sufficient but not necessary for studying discounted CTMDPs. For instance, we believe that the conditions imposed in \cite{Yan:2008}, which are different from the conditions imposed here and still allow unbounded transition rates and cost rates, could be also sufficient for us to obtain the standard results as presented in this paper. On the other hand, there exists research on CTMDPs (see \cite{Hu:2003}), whose study is only based on necessary conditions, which just requires that the underlying models are well defined (no explosion happens), and so are the expected total discounted costs (can be positive or negative infinity). In such a general setup, the authors of \cite{Hu:2003} obtain some nonstandard results for discounted CTMDPs in countable state and action spaces.

\section*{Appendix}
In this appendix, we establish some lemmas, and prove the main statements.

\begin{lemma}\label{MDPlem1}
Let a signed kernel $f(dy|x,t)$ on ${\cal B}(S)$ given $(x,t)\in S\times \mathbb{R}_+^0$ be fixed, and assume that it satisfies that following:
$f(\Gamma_S|x,t)\ge 0$ if
$\Gamma_S\in {\cal{B}}(S)$ and
$x\notin\Gamma_S$,  $f(S\setminus\{x\}|x,t)<\infty,$ and $f(S|x,t)=0$. Here, we put $F_x(t)\defi f(S\setminus\{x\}|x,t)<\infty$. Suppose there exist constants $\rho \ne 0,$ $b\ge 0$ and a measurable function $w(x)\ge0$ on $S$ such that $\int_Sf(dy|x,t)w(y)\le \rho w(x)+b,\forall~x\in S.$
Then
\begin{eqnarray*}
h(s,x,t)\ge\int_{s}^{t}\int_{S\setminus\{x\}}e^{-\int_{s}^{u}F_x(v)dv}f(dy|x,u)h(u,y,t)du
+e^{-\int_{s}^tF_x(v)dv}w(x),
\end{eqnarray*}
where $h$ is a nonnegative function defined by
\begin{eqnarray}\label{MDPlem1function}
h(s,x,t)\defi e^{\rho(t-s)}w(x)+\frac{b}{\rho}(e^{\rho(t-s)}-1),  \forall~0\le s\le t,x\in S.
\end{eqnarray}
\end{lemma}

\noindent{\textbf{Proof:}}
Straightforward calculations result in
\begin{eqnarray*}
&&\int_{s}^t\left\{e^{-\int_s^uF_x(v)dv}\int_{S\setminus\{x\}}f(dy|x,u)h(u,y,t)\right\}du+e^{-\int_s^tF_x(v)dv}w(x)\\
&=&\int_s^te^{-\int_s^uF_x(v)dv}e^{\rho(t-u)}\left(\int_{S}f(dy|x,u)w(y)-f(\{x\}|x,u)w(x)\right)du\\
&&+\frac{b}{\rho}\int_s^te^{-\int_s^uF_x(v)dv}e^{\rho(t-u)}F_x(u)du\\
&&-\frac{b}{\rho}\int_s^te^{-\int_s^uF_x(v)dv}F_x(u)du+e^{-\int_s^tF_x(v)dv}w(x)\\
&\le&\int_s^te^{-\int_s^uF_x(v)dv}e^{\rho(t-u)}\left(\rho w(x)+b+F_x(u)w(x)\right)du\\
&&+\frac{b}{\rho}\int_s^te^{-\int_s^uF_x(v)dv}e^{\rho(t-u)}F_x(u)du\\
&&-\frac{b}{\rho}\int_s^te^{-\int_s^uF_x(v)dv}F_x(u)du+e^{-\int_s^tF_x(v)dv}w(x).\\
\end{eqnarray*}
The rest of this proof now becomes identical to the one
of \cite[Lem.3.2(a), p.239]{Guo:2003IEEE}.
\hfill $\Box$
\bigskip

\begin{corollary}
Suppose Condition \ref{MDPRegularitycons}(b) is satisfied. If $\rho$ coming from Condition \ref{MDPRegularitycons} is strictly positive, then
\begin{eqnarray}\label{MDPneinproofLem2}
h(s,x,\tilde{t})&=&h(0,x,\tilde{t}-s)\nonumber\\
&\ge&\int_{s}^{\tilde{t}}\left\{e^{-\int_s^u\Lambda^l(S|x_0,\theta_1,\dots,\theta_l,x,v)dv}\int_{S}\Lambda^l(dy|x_0,\theta_1,\dots,\theta_l,x,u)h(u,y,\tilde{t})\right\}du\nonumber\\
&&+e^{-\int_s^{\tilde{t}}\Lambda^l(S|x_0,\theta_1,\dots,\theta_l,x,v)dv}w(x), \forall~x\in S,0\le s\le \tilde{t}<\infty,l\in\mathbb{Z}_+^0,
\end{eqnarray}
where $h$ is given in (\ref{MDPlem1function}).
\end{corollary}

\noindent{\textbf{Proof:}}
Let $l\in\mathbb{Z}_+^0$ be arbitrarily fixed. Consider the signed kernel on ${\cal{B}}(S)$
given $(x,u)\in S\times \mathbb{R}_+^0,$ defined by $\forall~\Gamma_S\in{\cal B}(S),$
\begin{eqnarray*}
g_l(\Gamma_S|x,u)\defi\left\{\begin{array}{ll}\displaystyle
\Lambda^l(\Gamma_S|x_0,\theta_1,\ldots,\theta_l,x,u) & \mbox{ if } x\notin \Gamma_S; \\ \\
\displaystyle -\Lambda^l(S|x_0,\theta_1,\ldots,\theta_l,x,u) &
\mbox{ if } \Gamma_S=\{x\},
\end{array}\right.
\end{eqnarray*}
where $\Lambda^l$ is defined in
(\ref{MDPexplicitintensity}). It can be easily verified that all the conditions in Lemma \ref{MDPlem1} are satisfied by $b\ge0,\rho>0,w(x)\ge1$ (coming from Condition \ref{MDPRegularitycons}) and this signed kernel $g_l(\cdot|x,u).$ Now the statement follows from Lemma \ref{MDPlem1}.
\hfill $\Box$
\bigskip

\begin{lemma}\label{MDPlem2}
Suppose Condition \ref{MDPRegularitycons}(b) is satisfied. Then under any policy
$\pi,$ $\forall~x\in S,m=0,1,2,\dots,$
\begin{eqnarray*}E_x^{\pi}[w(\xi_{t})I\{t<T_{m+1}\}]\le (e^{\rho
t}w(x)+\frac{b}{\rho}(e^{\rho
t}-1))I\{\rho>0\}+(w(x)+bt)I\{\rho=0\}.
\end{eqnarray*}
Here, constants $b,\rho$ and function $w$ come from Condition \ref{MDPRegularitycons}\footnote{In this lemma, we temporarily ignore Remark \ref{MDPrho}.}.
\end{lemma}

\noindent{\textbf{Proof:}}
Suppose $\rho>0.$ As for the statement, we prove the following slightly stronger result\footnote{Throughout this proof, this result is referred to as the ``stronger statement''.}, i.e., $\forall~m\in \mathbb{Z}_+^0,x\in  S,n=0,1,\dots,m,$
\begin{eqnarray*}
E_{x}^\pi\left[w(\xi_t)I\{t<T_{m+1}\}|{\cal{F}}_{T_{m-n}}\right]&\le&
I\{T_{m-n}\le
t\}h(T_{m-n},x_{m-n},t)\\
&&+\sum_{k=1}^{m-n}I\{T_{k-1}\le t<T_k\} w(x_{k-1}),
\end{eqnarray*}
where ${\cal{F}}_{T_{m-n}}\defi \sigma(x_i,T_i:i\in\mathbb{Z}_+^0, 0\le i\le m-n).$

This stronger statement is proved inductively.

Consider $n=0.$
On the set
$\{T_m\le t\},$ equation (\ref{MDPconstructprob}) implies
\begin{eqnarray}\label{MDPadditionaleqn}
P_x^{\pi}(\theta_{m+1}>t-T_m|{\cal{F}}_{T_m})
&=&e^{-\int_{0}^{t-T_m}\Lambda^m(S|h_m,v)dv}.
\end{eqnarray}
By the properties of conditional expectations and (\ref{MDPadditionaleqn}), we have
\begin{eqnarray*}
E_{x}^\pi\left[w(\xi_{t})I\{t<T_{m+1}\}\right|{\cal{F}}_{T_m}]&=&E_{x}^\pi\left[(I\{T_m\le
t\}+I\{T_m>t\})w(\xi_{t}) I\{t<T_{m+1}\}\right|{\cal{F}}_{T_m}]\\
&=&I\{T_{m}\le
t\}w(x_m) P_{x}^\pi(\theta_{m+1}>t-T_m|{\cal{F}}_{T_m})\\
&&+\sum_{k=1}^mI\{T_{k-1}\le t<T_k\}w(x_{k-1})\\
&=&I\{T_{m}\le
t\}w(x_m)e^{-\int_{0}^{t-T_m}\Lambda^m(S|h_m,v)dv}\\
&&+\sum_{k=1}^mI\{T_{k-1}\le t<T_k\}w(x_{k-1})\\
&\le&
I\{T_m\le t\}h(T_m,x_m,t)+\sum_{k=1}^mI\{T_{k-1}\le t<T_k\}w(x_{k-1}),
\end{eqnarray*}
where the last inequality follows from (\ref{MDPneinproofLem2}).

Now suppose the stronger statement
holds,  $\forall~0\le n<m$.

Consider the case of $n+1.$ By the properties of conditional expectations, the inductive supposition
and (\ref{MDPadditionaleqn}), we have
\begin{eqnarray*}
&&E_x^\pi\left[w(\xi_t)I\{t<T_{m+1}\}|{\cal{F}}_{T_{m-n-1}}\right]=E_x^\pi\left[E_x^\pi\left[w(\xi_t)I\{t<T_{m+1}\}|{\cal{F}}_{T_{m-n}}\right]|{\cal{F}}_{T_{m-n-1}}\right]\\
&\le& E_x^\pi\left[I\{T_{m-n}\le t\}h(T_{m-n},x_{m-n},t)
+\sum_{k=1}^{m-n}I\{T_{k-1}\le t<T_k\} w(x_{k-1})|{\cal{F}}_{T_{m-n-1}}\right]\\
&=&E_{x}^\pi\left[I\{T_{m-n-1}\le t\}I\{T_{m-n}\le
t\}h(T_{m-n},x_{m-n},t)|{\cal{F}}_{T_{m-n-1}}\right]\\
&&+E_{x}^\pi\left[I\{T_{m-n-1}\le
t\}\sum_{k=1}^{m-n}I\{T_{k-1}\le
t<T_{k}\}w(x_{k-1})|{\cal{F}}_{T_{m-n-1}}\right]\\
&&+E_{x}^\pi\left[I\{T_{m-n-1}>t\}\sum_{k=1}^{m-n}I\{T_{k-1}\le
t<T_k\}w(x_{k-1})|{\cal{F}}_{T_{m-n-1}}\right]\\
&=&I\{T_{m-n-1}\le t\}\left\{ \int_0^{t-T_{m-n-1}}\left\{\vphantom{\int_0^{t-T_{m-n-1}}}
 e^{-\int_{0}^u\Lambda^{m-n-1}(S|h_{m-n-1},v)dv}\right.\right.\\
&&\left.\left.\times\int_{S\setminus\{x_{m-n-1}\}}\Lambda^{m-n-1}(dy|h_{m-n-1},u)h(T_{m-n-1}+u,y,t)\right\}du\right.\\
&&\left.\vphantom{\int_0^{t-T_{m-n-1}}}+e^{-\int_{0}^{t-T_{m-n-1}}\Lambda^{m-n-1}(S|h_{m-n-1},v)dv}w(x_{m-n-1})\right\}+\sum_{k=1}^{m-n-1}I\{T_{k-1}\le
t<T_k\}w(x_{k-1})\\
&\le&
I\{T_{m-n-1}\le
t\}h(T_{m-n-1},x_{m-n-1},t)+\sum_{k=1}^{m-n-1}I\{T_{k-1}\le
t<T_k\} w(x_{k-1}),
\end{eqnarray*}
where the last inequality follows from (\ref{MDPneinproofLem2}).

Hence, the stronger statement holds. It remains to put
$n=m$ in the stronger statement to obtain Lemma \ref{MDPlem2}
for the case of $\rho>0$.

The statement corresponding to the case of $\rho=0$ follows from the fact of
$\lim_{\hat{\rho}\downarrow 0}\{e^{\hat{\rho} t}w(x)+\frac{b}{\hat{\rho}}(e^{\hat{\rho}
t}-1)\}=w(x)+bt.$ Here, we emphasize that if Condition \ref{MDPRegularitycons} is satisfied by $\rho=0,$ it is also satisfied by any arbitrarily fixed $\hat{\rho}>0.$
\hfill $\Box$
\bigskip

\begin{lemma}\label{MDPextralemma}
Suppose Condition \ref{MDPRegularitycons} is satisfied. For any fixed $l\in\mathbb{Z}_+^0,$ consider the modified transition rates defined by
$$\tilde q_l(\cdot|x,a)\defi\left\{\begin{array}{ll}
q(\cdot|x,a), & \mbox{ if } x\in S_l; \\
0, & \mbox{ if } x\in S\setminus S_l. \end{array}\right. $$ Their
corresponding probabilities and expectations are denoted by $P_x^{\pi,l}$ and $E_x^{\pi,l}.$ Then under any policy $\pi,$ $\forall~x\in S, t\ge 0,$
\begin{eqnarray}\label{MDPauxi}
\lim_{l\rightarrow\infty}\tilde{P}_x^{\pi,l}(\xi_t\in S\setminus S_l)=0,
\end{eqnarray}
where $S_l$ is defined in Condition \ref{MDPRegularitycons}(a).
\end{lemma}

\noindent{\textbf{Proof:}}
Throughout this proof, let $x\in S$ and $t\ge 0$ be arbitrarily fixed. Under Condition \ref{MDPRegularitycons}, we have that $\forall~\epsilon>0,\exists~J(\epsilon)>0: \forall~l\ge J(\epsilon),$
\begin{equation}\label{MDPproofofthm1eqn2}
\inf_{y\in S\setminus
S_{l}}w(y)>\frac{e^{\tilde{\rho}t}w(x)+\frac{b}{\tilde{\rho}}(e^{\tilde{\rho}t}-1)}{\epsilon},
\end{equation}
where $\tilde{\rho}\defi \rho +1.$

Suppose the statement of this lemma does
not hold, i.e., $\exists~\epsilon>0:\forall~ L>0,\exists~
l\ge \max\{L,~J(\epsilon)\}:$
\begin{equation}\label{MDPproofofthm1opposite}
\tilde{P}_x^{\pi,l}(\xi_t\in S\setminus S_l)>\epsilon.
\end{equation}
At the same time, necessarily, (\ref{MDPproofofthm1eqn2}) holds as well. On the one hand, by using Lemma \ref{MDPlem2}\footnote{If Condition \ref{MDPRegularitycons} is satisfied by $\rho$
and $q$, then it is also satisfied by $\tilde{\rho}$ and $\tilde{q}$, where we recall $\tilde{\rho}=1+\rho.$} and the fact of $\sup_{x\in S}\sup_{a\in
A(x)}\tilde{q}_{x}(a)\le\sup_{x\in {S_l}}\bar{q}_x<\infty$ (see
Condition \ref{MDPRegularitycons}), we have
\begin{eqnarray}\label{MDPproofofthmcontra}
\tilde{E}^{\pi,l}_x\left[w(\xi_t)\right]&=&\tilde{E}^{\pi,l}_x\left[w(\xi_t)\sum_{m=0}^{\infty}I\{T_m\le
t<T_{m+1}\}\right]=\lim_{m\rightarrow\infty}\tilde{E}_x^{\pi,l}[w(\xi_t)I\{t<T_{m+1}\}]
\nonumber\\
&\le&
e^{\tilde{\rho}t}w(x)+\frac{b}{\tilde{\rho}}(e^{\tilde{\rho}t}-1).
\end{eqnarray}
On the other hand, we have
\begin{eqnarray*}
\tilde{E}^{\pi,l}_x\left[w(\xi_t)\right]&=&\tilde{E}^{\pi,l}_x\left[w(\xi_t)|\xi_t\in S\setminus S_l\right]\tilde{P}^{\pi,l}_x(\xi_t\in S\setminus S_l)+\tilde{E}^{\pi,l}_x\left[w(\xi_t)|\xi_t\in S_l\right]\tilde{P}^{\pi,l}_x(\xi_t\in S_l)\\
&>&\inf_{y\in {S\setminus S_l}}w(y)\epsilon>
e^{\tilde{\rho}t}w(x)+\frac{b}{\tilde{\rho}}(e^{\tilde{\rho}t}-1),
\end{eqnarray*}
where the first inequality follows from ignoring the second term
in the first line and estimating the first term from below using
(\ref{MDPproofofthm1opposite}), and the last inequality is a result of
(\ref{MDPproofofthm1eqn2}). However, this contradicts
(\ref{MDPproofofthmcontra}). \hfill $\Box$
\bigskip

\noindent{\textbf{Proof of Theorem \ref{MDPregularitythm}:}}
(a) From (\ref{MDPconstructprob}), we clearly have that $\forall~l\in \mathbb{Z}_+^0,t\ge0$,
\begin{eqnarray}\label{MDPnewextra}
&&P_x^\pi\left((\xi_t=x_\infty)\bigcup((\xi_t\ne x_\infty)\bigcap (\mbox{the process visits $S\setminus S_{l}$ at least once on $[0,t]$}))\right)\nonumber\\
&=&1-P_x^\pi(\forall~\tilde{t}\in[0,t],\xi_{\tilde{t}}\in S_l)=1-\tilde{P}_x^{\pi,l}(\xi_t\in S_l)\nonumber\\
&=& \tilde{P}_x^{\pi,l}\left((\xi_t=x_\infty)\bigcup (\xi_t\in S\setminus S_l)\right)= \tilde{P}_x^{\pi,l}\left(\xi_t\in S\setminus S_l\right).
\end{eqnarray}
Here, we have repeatedly used the fact of $\sup_{x\in S}\sup_{a\in
A(x)}\tilde{q}_{x}(a)\le\sup_{x\in {S_l}}\bar{q}_x<\infty$, so that $\tilde{P}_x^{\pi,l}(T_\infty=\infty)=1$.
By using Lemma \ref{MDPextralemma}, (\ref{MDPnewextra}) and the fact that $(S\setminus S_l)_{l\in \mathbb{Z}_+^0}$ is a decreasing system, we have $\forall~t\ge 0,$
\begin{eqnarray*}
P_x^\pi\left(\forall~l\in\mathbb{Z}_+^0,(\xi_t=x_\infty)\bigcup((\xi_t\ne x_\infty)\bigcap (\mbox{the process visits $S\setminus S_{l}$ at least once on $[0,t]$}))\right)=0,
\end{eqnarray*}
which is equivalent to
\begin{eqnarray*}
P_x^\pi\left(\exists~l\in\mathbb{Z}_+^0,(\xi_t\ne x_\infty)\bigcap((\xi_t= x_\infty)\bigcup(\forall~\tilde{t}\in[0,t],\xi_{\tilde{t}}\in S_l))\right)=1,
\end{eqnarray*}
i.e., for each $t\ge 0$, $P_x^\pi(\exists~l\in\mathbb{Z}_+^0, \forall~\tilde{t}\in[0,t],\xi_{\tilde{t}}\in S_l)=1$. However, if
$\xi_{\tilde{t}}\in S_l$ on $[0,t]$ a.s., then $T_\infty>t,$ a.s., i.e, $P_x^\pi(T_\infty>t)=1.$ Since $t\ge 0$ is arbitrary,
this leads to $P_x^\pi(T_\infty=\infty)=1$ and $P_x^\pi(\xi_t\in S)=1, \forall~t\ge 0.$ The statement regarding $E_x^\pi[w(\xi_t)]$ follows from this, Lemma \ref{MDPlem2} and that $\forall~t\ge 0,$
\begin{eqnarray*}
E_x^\pi\left[w(\xi_t)\right]=E_x^\pi\left[w(\xi_t)\sum_{m=0}^\infty
I\{T_m\le t<T_{m+1}\}\right]=\lim_{m\rightarrow\infty}E_x^\pi\left[w(\xi_t)I\{t<T_{m+1}\}\right].
\end{eqnarray*}

(b) By definition, we have
$
V_{0}(x,\pi)\defi E_{x}^\pi\left[\int_0^\infty e^{-\alpha
t}\int_{A}c_0(\xi_{t-},a)\pi(da|\omega,t)dt\right]
$. Then, using Condition
\ref{MDPfinitenesscons}(b,c) and Theorem
\ref{MDPregularitythm}(a), we obtain
\begin{eqnarray*}
V_0(x,\pi)&\ge& -E_x^\pi\left[\int_0^\infty e^{-\alpha
t}(Mw(\xi_{t})+c)dt\right]=-\int_{0}^\infty e^{-\alpha t}(ME_x^\pi[w(\xi_{t})]+c)dt \\
 &\ge& -\int_0^\infty e^{-\alpha t}(M(e^{\rho
t}w(x)+\frac{b}{\rho}(e^{\rho t}-1))+c)dt=-\frac{M(\alpha w(x)+b)}{\alpha(\alpha-\rho)}-\frac{c}{\alpha}.
\end{eqnarray*}
With Condition \ref{MDPfinitenesscons}(a) in mind, the statement for $V_0(\pi)=\int_S V_0(x,\pi)\gamma(dx)$
follows.
\hfill $\Box$
\bigskip

\noindent\textbf{Proof of Theorem \ref{MDPkolmogorov}:} (a) Similarly to $\mu$ and $\nu$ (defined by (\ref{MDPrandommeasure1}) and (\ref{MDPrandommeasure2})), let us define the
following two random measures :
\begin{eqnarray*}
\tilde{\mu}(\omega,dt,\Gamma)\defi \sum_{m\ge
1}I\{T_{m}<\infty\}I\{x_{m-1}\in \Gamma\}I\{T_{m}\in dt\}, \forall~\Gamma\in {\cal B}(S)
\end{eqnarray*}
and
\begin{eqnarray*}
\tilde{\nu}(\omega,dt,\Gamma)\defi
\int_{A}\pi(da|\omega,t)q(S\setminus\{\xi_{t-}\}|\xi_{t-},a)I\{\xi_{t-}\in
\Gamma\}dt, \forall~\Gamma\in {\cal{B}}(S).
\end{eqnarray*}
It is shown in the proof of \cite[Lem.4]{Kitaev:1986} that
$\tilde{\nu}$ is the dual predictable projection of $\tilde{\mu},$ i.e.,
for any nonnegative ${\cal{P}}\times {\cal{B}}(S)$\footnote{Here, we clarify that ${\cal{P}}\times {\cal{B}}(S)$ denotes the product
$\sigma$-algebra, rather than the Cartesian product.}-measurable function $Y(\omega,t,x)$,
\begin{eqnarray*}
E_{x}^\pi\left[\int_{0}^\infty\int_{S}\tilde{\mu}(dt,dy)Y(t,y)\right]=E_{x}^\pi\left[\int_{0}^\infty\int_{S}\tilde{\nu}(dt,dy)Y(t,y)\right],
\end{eqnarray*}
see \cite[Chap.4, Sec.5]{Kitaev:1995} for more details.  Now it immediately follows that
$E_x^\pi\left[\tilde{\mu}((0,t],\Gamma)\right]<\infty,$ because by using Condition \ref{MDPRegularitycons}(c) and the definition of $\Gamma$ given in the statement of this theorem, we have
\begin{eqnarray}\label{MDPfinitekitaev}
E_x^\pi\left[\tilde{\nu}((0,t],\Gamma)\right]&=&E_x^\pi\left[\int_0^t\int_A\pi(da|\omega,u)q_{\xi_{u-}}(a)I\{\xi_{u-}\in\Gamma\}du\right]\nonumber\\
&\le&
t\sup_{y\in S_l}\bar{q}_{y}<\infty.
\end{eqnarray}
On the other hand, by Theorem \ref{MDPregularitythm}, $\mu((0,t],\Gamma)$
and $\tilde{\mu}((0,t],\Gamma)$ are a.s. finite. Then it
follows from their definitions that
$|\mu((0,t],\Gamma)-\tilde{\mu}((0,t],\Gamma)|\le 1$ a.s.. Therefore, $E_x^\pi\left[\mu((0,t],\Gamma)\right]<\infty$.
Consequently, it is legal to take expectations
in the both sides of the following obviously valid equation
$$
I\{\xi_t\in \Gamma\}=I\{\xi_{0}\in
\Gamma\}+\mu((0,t],\Gamma)-\tilde{\mu}((0,t],\Gamma) \mbox{ a.s.},
$$
from which the statement follows.

(b) The reasoning for proving part (a) of this theorem can be repeated, except that now one needs replace the argument for
(\ref{MDPfinitekitaev}) by the following:
\begin{eqnarray*}
E_x^\pi\left[\tilde{\nu}((0,t],\Gamma)\right]&=&E_x^\pi\left[\int_0^t\int_A\pi(da|\omega,u)q_{\xi_{u-}}(a)I\{\xi_{u-}\in\Gamma\}du\right]\nonumber\\
&\le&
E_x^\pi\left[\int_0^tLw(\xi_{u-})I\{\xi_{u-}\in\Gamma\}du\right]\nonumber\\
&\le& L\int_{0}^tE_x^\pi\left[w(\xi_{u})\right]du<\infty,
\end{eqnarray*}
where the second inequality follows from Condition \ref{MDPstrongregularitycons}, and the last inequality is due to Theorem \ref{MDPregularitythm}.
\hfill $\Box$
\bigskip

\noindent \textbf{Proof of Theorem \ref{MDPDynkinsformula}:}
Step 1. We prove
that equation (\ref{MDPDynkinformula1}) holds for $r(x)\defi u(x)I\{x\in S_l\}$, where $S_l$ is defined in Condition \ref{MDPRegularitycons}.

We obviously have
\begin{eqnarray}\label{eqm1}
&&\int_Sw'(y)E_x^\pi\left[\int_0^t\int_A\pi(da|\omega,v)q(dy\setminus
\{\xi_v\}|\xi_v,a)dv\right]\nonumber\\
&=&E_x^\pi\left[\int_0^t\int_A\pi(da|\omega,v)\int_Sw'(y)q(dy\setminus
\{\xi_v\}|\xi_v,a)dv\right]\\
&=&E_x^\pi\left[\int_0^t\int_A\pi(da|\omega,v)\int_Sw'(y)\left\{q(dy|\xi_v,a)-q(\{\xi_v\}|\xi_v,a)I\{\xi_v\in
dy\}\right\}dv\right]<\infty.\nonumber
\end{eqnarray}
Indeed, by Condition \ref{MDPdpcons}(a,b) and Theorem \ref{MDPregularitythm}(a),
\begin{eqnarray*}
&&E_x^\pi\left[\int_0^t\int_A\pi(da|\omega,v)\int_Sw'(y)q(dy|\xi_v,a)dv\right]\le E_x^\pi\left[\int_0^t\int_A\pi(da|\omega,v)(\rho'w'(\xi_v)+b')dv\right]\\
&\le& L'\rho'\int_0^tE_x^\pi\left[ w(\xi_v)\right]dv+b't<\infty,
\end{eqnarray*}
and
\begin{eqnarray}\label{eqm2}
&&E_x^\pi\left[\int_0^t\int_A\pi(da|\omega,v)w'(\xi_v)|q(\{\xi_v\}|\xi_v,a)|dv\right] =E_x^\pi\left[\int_0^t\int_A\pi(da|\omega,v)w'(\xi_v)q_{\xi_v}(a)dv\right]\nonumber\\
&\le& L'\int_0^tE_x^\pi\left[w(\xi_v)\right]dv<\infty.
\end{eqnarray}

It follows from the previous calculations that
\begin{eqnarray*}
&&\int_Sr(y)E_x^\pi\left[\int_0^t\int_A\pi(da|\omega,v)q(dy\setminus\{\xi_v\}|\xi_v,a)dv\right]\\
&\le& ||r||_{w'}\int_Sw'(y)E_x^\pi\left[\int_0^t\int_A\pi(da|\omega,v)q(dy\setminus
\{\xi_v\}|\xi_v,a)dv\right]<\infty,
\end{eqnarray*}
and
\begin{eqnarray*}\label{MDPdynkinfinite2}
E_x^\pi\left[\int_0^t\int_A\pi(da|\omega,v)q_{\xi_v}(a)r(\xi_v)dv\right]<\infty.
\end{eqnarray*}
Now in order to establish equation (\ref{MDPDynkinformula1}) for $r(x)=u(x)I\{x\in S_l\}$, one only needs integrate $r(x)$ over $S$ with respect to $P_x^\pi(\xi_t\in \cdot)$ and use Theorem \ref{MDPkolmogorov}.

Step 2. We prove that equation (\ref{MDPDynkinformula1}) holds for any $u(x)\in \textbf{B}_{w'}(S)$. By putting $S_{-1}\defi
\emptyset$ and observing $E_x^\pi\left[\sum_{l=-1}^\infty
|u(\xi_t)|I\{\xi_t\in S_{l+1}\setminus S_l\}\right]<\infty,$ we have
\begin{eqnarray*}
&&E_x^\pi\left[u(\xi_t)\right]-u(x)=E_x^\pi\left[\sum_{l=-1}^\infty u(\xi_t)I\{\xi_t\in
S_{l+1}\setminus S_l\}\right]-\sum_{l=-1}^\infty u(x)I\{x\in
S_{l+1}\setminus S_l\}
\\
&=&\sum_{l=-1}^{\infty}E_x^\pi\left[u(\xi_t)I\{\xi_t\in
S_{l+1}\setminus S_l\}\right]-\sum_{l=-1}^\infty u(x)I\{x\in
S_{l+1}\setminus S_l\}\\
&=&\sum_{l=-1}^{\infty}\left\{E_x^\pi\left[u(\xi_t)I\{\xi_t\in
S_{l+1}\setminus S_l\}\right]- u(x)I\{x\in S_{l+1}\setminus
S_l\}\right\}\\
&=&\sum_{l=-1}^\infty
\left\{E_x^\pi\left[\int_0^t\int_S\int_A\pi(da|\omega,v)q(dy|\xi_v,a)u(y)I\{y\in
S_{l+1}\setminus S_l\}\right]\right\}\\
&=&E_x^\pi\left[\int_0^t\int_S\int_A\pi(da|\omega,v)q(dy|\xi_{v},a)u(y)dv\right],
\end{eqnarray*}
where the second last equality follows from formally applying the result obtained in Step 1 of this proof, i.e., (\ref{MDPDynkinformula1}) holds for
$r(x).$ The involved interchange of the order of integrations, summations and expectations is legal, as can be easily verified similarly to (\ref{eqm1}) and (\ref{eqm2}).

Step 3. We prove that equation (\ref{MDPDynkinformula2}) holds for any $u(x)\in \textbf{B}_{w'}(S).$
In this proof, we repeatedly apply (\ref{MDPDynkinformula1}) to $E_x^\pi[u(\xi_t)].$ On the one hand, we have
\begin{eqnarray*}
\mbox{LHS of (\ref{MDPDynkinformula2})}&=&e^{-\alpha
t}\left\{u(x)+E_x^\pi\left[\int_0^t\int_S\int_A\pi(da|\omega,v)q(dy|\xi_v,a)u(y)dv\right]\right\}-u(x)\\
&=&e^{-\alpha
t}E_x^\pi\left[\int_0^t\int_S\int_A\pi(da|\omega,v)q(dy|\xi_v,a)u(y)dv\right]+u(x)(e^{-\alpha
t}-1).
\end{eqnarray*}
On the other hand, we have the following two observations. Firstly,
\begin{eqnarray*}
&&E_x^\pi\left[\int_0^te^{-\alpha v}(-\alpha
u(\xi_v))dv\right]=-\alpha \int_0^t e^{-\alpha
v}E_x^\pi\left[u(\xi_v)\right]dv\\
&=&-\alpha\int_0^te^{-\alpha
v}\left\{u(x)+E_x^\pi\left[\int_0^v\int_S\int_A\pi(da|\omega,r)q(dy|\xi_r,a)u(y)dr\right]\right\}dv\\
&=&(e^{-\alpha t}-1)u(x)-\alpha\int_0^te^{-\alpha
v}E_x^\pi\left[\int_0^v\int_S\int_A\pi(da|\omega,r)q(dy|\xi_r,a)u(y)dr\right]dv\\
&=&(e^{-\alpha t}-1)u(x)-\alpha
E_x^\pi\left[\int_0^t\left\{e^{-\alpha
v}\int_0^v\int_S\int_A\pi(da|\omega,r)q(dy|\xi_r,a)u(y)dr\right\}dv\right]
\end{eqnarray*}
where the interchange of the order of integrals in the first and the last
equalities is legal, because evidently, $\forall~u\in
\textbf{B}_{w'}(S),$
$
E_x^\pi\left[\int_0^te^{-\alpha v}\alpha |u(\xi_v)|dv\right]<\infty
$
and
$$\int_0^te^{-\alpha
v}E_x^\pi\left[\int_0^v\int_S\int_A\pi(da|\omega,r)q(dy|\xi_r,a)|u|(y)dr\right]dv<\infty.$$ Secondly, integration by parts results in
\begin{eqnarray*}
&&E_x^\pi\left[\int_0^te^{-\alpha
v}\int_S\int_A\pi(da|\omega,v)q(dy|\xi_v,a)u(y)dv\right]\\
&=&E_x^\pi\left[e^{-\alpha
t}\int_0^t\int_S\int_A\pi(da|\omega,r)q(dy|\xi_r,a)u(y)dr\right]\\
&&+\alpha E_x^\pi\left[\int_0^t e^{-\alpha
v}\int_0^v\int_S\int_A\pi(da|\omega,r)q(dy|\xi_r,a)u(y)dr~dv\right].
\end{eqnarray*}
These two observations, together with the expression for LHS of (\ref{MDPDynkinformula2}) obtained in the above, finally lead to
\begin{eqnarray*}
&&\mbox{RHS of
(\ref{MDPDynkinformula2})}\\
&=&E_x^\pi\left[\int_0^te^{-\alpha
v}(-\alpha u(\xi_v))dv\right]+E_x^\pi\left[\int_0^te^{-\alpha
v}\int_S\int_A\pi(da|\omega,v)q(dy|\xi_v,a)u(y)dv\right]\\
&=&(e^{-\alpha t}-1)u(x)+E_x^\pi\left[e^{-\alpha
t}\int_0^t\int_S\int_A\pi(da|\omega,r)q(dy|\xi_r,a)u(y)dr\right]
=\mbox{ LHS of (\ref{MDPDynkinformula2})},
\end{eqnarray*}
as required.
\hfill $\Box$
\bigskip

\begin{lemma}\label{MDPmeaslem}
Suppose Condition \ref{MDPRegularitycons}(b) and Condition \ref{MDPmeas2} are satisfied. Then
$\forall~u\in \textbf{B}_{w}(S),$  function $v$ given by
$$
v(x)\defi \inf_{a\in
A(x)}\left\{\frac{c_0(x,a)}{\alpha+1+\bar{q}_x}+\frac{1+\bar{q}_x}{\alpha+1+\bar{q}_x}\int_S
u(y)\left(\frac{q(dy|x,a)}{1+\bar{q}_x}+I\{x\in
dy\}\right)\right\}
$$
is measurable in $x\in S$.
\end{lemma}

\par\noindent\textbf{Proof:}
By Remark \ref{MDPremarkextra}, Condition \ref{MDPRegularitycons}(b) and Condition \ref{MDPmeas2},  we refer to \cite[Lem.8.3.7(a)]{Hernandez-Lerma:1999} for that $\forall~u\in \textbf{B}_w(S),x\in S,$ function\footnote{It can be easily verified that $\forall~(x,a)\in K,$ $\left(\frac{q(dy|x,a)}{1+\bar{q}_x}+I\{x\in dy\}\right)$ is a probability measure on $(S,{\cal B}(S))$.} $\int_S u(y)\left(\frac{q(dy|x,a)}{1+\bar{q}_x}+I\{x\in dy\}\right)$ is
continuous in $a\in A(x).$  It follows from this and Condition \ref{MDPmeas2}(c) that $\forall~ x\in S, u\in \textbf{B}_{w}(S),$ function
\begin{eqnarray*}
\frac{c_0(x,a)}{\alpha+1+\bar{q}_x}+\frac{1+\bar{q}_x}{\alpha+1+\bar{q}_x}\int_S
u(y)\left(\frac{q(dy|x,a)}{1+\bar{q}_x}+I\{x\in
dy\}\right)
\end{eqnarray*}
 is lower semicontinuous in $a\in A(x).$
By \cite[Prop.7.29]{Bertsekas:1978}, $\forall~u\in \textbf{B}_{w}(S),$ function
\begin{eqnarray*}
\frac{c_0(x,a)}{\alpha+1+\bar{q}_x}+\frac{1+\bar{q}_x}{\alpha+1+\bar{q}_x}\int_S
u(y)\left(\frac{q(dy|x,a)}{1+\bar{q}_x}+I\{x\in
dy\}\right)
\end{eqnarray*}
is measurable\footnote{We emphasize that by Remark \ref{MDPremarkextra}, we have that $\bar{q}_x$ is measurable on $S.$} on $K$. Now it remains to apply
\cite[D.5 Prop.]{Hernandez-Lerma:1996} (see also
\cite[Prop.7.33]{Bertsekas:1978}) for the statement of this lemma.
\hfill $\Box$
\bigskip

\noindent\textbf{Proof of Theorem \ref{MDPguosprop}:}
Throughout this proof, $x\in S$ is arbitrarily fixed. Due to Lemma \ref{MDPmeaslem}, functions $u^{(n)},n=0,1,2,\dots$ are measurable. Now the proof goes in steps.

Step 1. We prove that $\{u^{(n)},n=0,1,\dots\}$ is a non-increasing sequence.

Straightforward calculations result in
\begin{eqnarray*}
&&u^{(1)}(x)=\inf_{a\in A(x)}\left\{\frac{c_0(x,a)}{\alpha+1+\bar{q}_x}+\frac{1+\bar{q}_x}{\alpha+1+\bar{q}_x} \int_S u^{(0)}(y)\left(\frac{q(dy|x,a)}{1+\bar{q}_x}+I\{x\in dy\}\right)\right\}\\
&=&\inf_{a\in A(x)}\left\{\frac{c_0(x,a)}{\alpha+1+\bar{q}_x}+\frac{1+\bar{q}_x}{\alpha+1+\bar{q}_x}\int_S\left(\frac{M(\alpha w(y)+b)}{\alpha(\alpha-\rho)}+\frac{c}{\alpha}\right)\left(\frac{q(dy|x,a)}{1+\bar{q}_x}+I\{x\in dy\}\right)\right\}\\
&\le&\inf_{a\in A(x)}\left\{\frac{c_0(x,a)}{\alpha+1+\bar{q}_x}\right\}\\
&&+\frac{1+\bar{q}_x}{\alpha+1+\bar{q}_x}\sup_{a\in A(x)}\left\{\int_S\left(\frac{M(\alpha w(y)+b)}{\alpha(\alpha-\rho)}+\frac{c}{\alpha}\right)\left(\frac{q(dy|x,a)}{1+\bar{q}_x}+I\{x\in dy\}\right) \right\}\\
&\le& \frac{Mw(x)+c}{\alpha+1+\bar{q}_x}+\frac{1+\bar{q}_x}{\alpha+1+\bar{q}_x}\left\{\frac{bM}{\alpha (\alpha-\rho)}+\frac{M(\rho w(x)+b)}{(\alpha-\rho)(1+\bar{q}_x)}+\frac{M w(x)}{\alpha-\rho}+\frac{c}{\alpha}\right\}=u^{(0)}(x),
\end{eqnarray*}
where the last inequality follows from Condition \ref{MDPRegularitycons}(b) and Condition \ref{MDPfinitenesscons}(c). Now the result of Step 1 follows from this and the monotonicity of the RHS of (\ref{eqm}) with respect to $u^{(n)}$.

Step 2. We prove that $\forall~n=0,1,\dots, |u^{(n)}(x)|\le \frac{M(\alpha w(y)+b)}{\alpha(\alpha-\rho)}+\frac{c}{\alpha}=u^{(0)}(x).$

On the one hand, the result of Step 1 implies that $\forall~n=0,1,\dots, u^{(n)}(x)\le  \frac{M(\alpha w(y)+b)}{\alpha(\alpha-\rho)}+\frac{c}{\alpha}.$ On the other hand, we have that
\begin{eqnarray*}
&&u^{(1)}(x)=\inf_{a\in A(x)}\left\{\frac{c_0(x,a)}{\alpha+1+\bar{q}_x}+\frac{1+\bar{q}_x}{\alpha+1+\bar{q}_x} \int_S u^{(0)}(y)\left(\frac{q(dy|x,a)}{1+\bar{q}_x}+I\{x\in dy\}\right)\right\}\\
&=& \inf_{a\in A(x)}\left\{\frac{c_0(x,a)}{\alpha+1+\bar{q}_x}+\frac{1+\bar{q}_x}{\alpha+1+\bar{q}_x}\int_S\left(\frac{M(\alpha w(y)+b)}{\alpha(\alpha-\rho)}+\frac{c}{\alpha}\right)\left(\frac{q(dy|x,a)}{1+\bar{q}_x}+I\{x\in dy\}\right)\right\}\\
&\ge& \inf_{a\in A(x)}\left\{\frac{c_0(x,a)}{\alpha+1+\bar{q}_x}\right\}\\
&&+\inf_{a\in A(x)} \left\{\frac{1+\bar{q}_x}{\alpha+1+\bar{q}_x}\int_S\left(\frac{M(\alpha w(y)+b)}{\alpha(\alpha-\rho)}+\frac{c}{\alpha}\right)\left(\frac{q(dy|x,a)}{1+\bar{q}_x}+I\{x\in dy\}\right)
\right\}\\
&\ge&-\frac{Mw(x)+c}{\alpha+1+\bar{q}_x}\\
&&+\frac{1+\bar{q}_x}{\alpha+1+\bar{q}_x}\inf_{a\in A(x)}\left\{\int_S-\left(\frac{M(\alpha w(y)+b)}{\alpha(\alpha-\rho)}+\frac{c}{\alpha}\right)\left(\frac{q(dy|x,a)}{1+\bar{q}_x}+I\{x\in dy\}\right)
\right\}\\
&=&-\frac{Mw(x)+c}{\alpha+1+\bar{q}_x}\\
&&-\frac{1+\bar{q}_x}{\alpha+1+\bar{q}_x}\sup_{a\in A(x)}\left\{\int_S\left(\frac{M(\alpha w(y)+b)}{\alpha(\alpha-\rho)}+\frac{c}{\alpha}\right)\left(\frac{q(dy|x,a)}{1+\bar{q}_x}+I\{x\in dy\}\right)
\right\}\\
&\ge& -\frac{Mw(x)+c}{\alpha+1+\bar{q}_x}-\frac{1+\bar{q}_x}{\alpha+1+\bar{q}_x}\left\{\frac{bM}{\alpha (\alpha-\rho)}+\frac{M(\rho w(x)+b)}{(\alpha-\rho)(1+\bar{q}_x)}+\frac{M w(x)}{\alpha-\rho}+\frac{c}{\alpha}\right\}=-u^{(0)}(x),
\end{eqnarray*}
where the second inequality is because of Condition \ref{MDPfinitenesscons}(c), $\frac{M(\alpha w(y)+b)}{\alpha(\alpha-\rho)}+\frac{c}{\alpha}\ge 0$ and the fact of $\frac{q(dy|x,a)}{1+\bar{q}_x}+I\{x\in dy\}$ being a probability measure, and the last inequality follows from Condition \ref{MDPRegularitycons}(b). This and an inductive argument lead to that $\forall~ n=0,1,\dots, u^{(n)}(x)\ge -\left(\frac{M(\alpha w(y)+b)}{\alpha(\alpha-\rho)}+\frac{c}{\alpha}\right).$ Thus, Step 2 is completed.

Now it follows from the results of Step 1 and Step 2 that $u^\ast(x)=\lim_{n\rightarrow\infty}u^{(n)}(x)$ exists and $u^{\ast}(x)\in \textbf{B}_{w}(S).$ The fact that $u^\ast$ solves the Bellman equation (\ref{MDPdpeqn}) can be verified in exactly the same way as in the proof of \cite[Lem.3.3(b)]{Guo:2007}, and its proof is thus omitted.
\hfill $\Box$
\bigskip

\begin{lemma}\label{MDPdplemexpressionforperformance}
Suppose Condition \ref{MDPRegularitycons},
Condition \ref{MDPfinitenesscons}(a,b), Condition \ref{MDPdpcons} and
Condition \ref{MDPmeas2} are
satisfied. Then under any policy $\pi,$
\begin{eqnarray}\label{MDPdplemeqnexpression}
V_0(\pi)&=&E_\gamma^{\pi}\left[\int_0^\infty
e^{-\alpha t}\int_A\pi(da|\omega,t)\left\{c_0(\xi_t,a)-\alpha
u(\xi_t)+\int_S q(dy|\xi_t,a)u(y)\right\}dt\right]\nonumber\\
&&+\int_S\gamma(dy)u(y),
\end{eqnarray}
where  $u\in \textbf{B}_{w'}(S)$ is an arbitrary function.
\end{lemma}

\noindent\textbf{Proof:}
By applying Dynkin's
formula (\ref{MDPDynkinformula2}) to $e^{-\alpha
t}E_\gamma^\pi\left[u(\xi_t)\right],$ we have
$$e^{-\alpha
t}E_\gamma^\pi\left[u(\xi_t)\right]=\int_S\gamma(dy)u(y)+E_\gamma^\pi\left[\int_0^te^{-\alpha
v}\int_A \pi(da|\omega,v)\left\{-\alpha
u(\xi_v)+\int_Sq(dy|\xi_v,a)u(y)\right\}dv\right].$$
The expectations of all particular summands are finite here. According to Theorem \ref{MDPregularitythm}(b) (see also its proof), we can formally
add $E_\gamma^\pi\left[
\int_0^te^{-\alpha v}\int_A\pi(da|\omega,v)c_0(\xi_v,a) dv\right]$
to the both sides of the above equation, and take the limit as $t\rightarrow \infty$.
We emphasize that $\lim_{t\rightarrow\infty}e^{-\alpha
t}E_\gamma^\pi\left[u(\xi_t)\right]=0$ because of Theorem \ref{MDPregularitythm}(a) and Condition \ref{MDPfinitenesscons}(b).
\hfill $\Box$
\bigskip

The next lemma can be established in exactly the same way as in the proof of \cite[Lem.5.3]{Guo:2007}.
\begin{lemma}\label{MDPDPinequality}
Suppose Condition \ref{MDPRegularitycons},
Condition \ref{MDPfinitenesscons}(a,b), Condition \ref{MDPdpcons} and
Condition \ref{MDPmeas2} are
satisfied. Then under any fixed Markov policy
$\pi,$ $\forall~x\in S$, the following assertions hold:
\par\noindent (a) If $u\in \textbf{B}_{w'}(S),$ and
$
\alpha u(x)\ge \int_A
\pi(da|x,t)c_0(x,a)+\int_S\int_A\pi(da|x,t)q(dy|x,a)u(y),\forall~x\in S,t\ge 0,$ then $u(x)\ge V_0(x,\pi).$
\par\noindent (b) If $u\in \textbf{B}_{w'}(S),
$
and
$
\alpha u(x)\le \int_A \pi(da|x,t)c_0(x,a)+\int_S\int_A\pi(da|x,t)q(dy|x,a)u(y),\forall~x\in S,t\ge 0,$ then $u(x)\le V_0(x,\pi).
$
\end{lemma}

\noindent\textbf{Proof of Theorem \ref{MDPbellmanthm}:} (a)
Using \cite[D.5 Prop.]{Hernandez-Lerma:1996} and the fact that $u^\ast$ solves the Bellman equation (\ref{MDPdpeqn}), we have that $\forall~\epsilon>0,\exists$
a deterministic stationary policy $\hat{\phi}:$
$$
c_0(x,\hat{\phi}(x))- \alpha
u^\ast(x)+\int_Sq(dy|x,\hat{\phi}(x))u^\ast(y)\le\alpha \epsilon, \forall~x\in S.
$$
It follows from this and Lemma \ref{MDPdplemexpressionforperformance} that
$V_0(\hat{\phi})\le\int_S\gamma(dy)u^\ast(y)+\epsilon,$ and thus\footnote{Here, we recall that $\epsilon>0$ is arbitrary.} $\inf_\phi V_0(\phi)\le
\int_S\gamma(dy)u^\ast(y).$ On the other hand, by Lemma \ref{MDPdplemexpressionforperformance}, we have that under any policy $\pi,$ $V_0(\pi)\ge\int_{S}\gamma(dy)u^\ast(y).$
Now it is evident that $\int_S\gamma(dy)u^\ast(y)=\inf_{\pi}V_0(\pi)=\inf_{\phi}V_0(\phi).$ The proof for the existence of a deterministic stationary optimal policy is identical (with few very minor modifications) to the one of \cite[Thm.3.3(c)]{Guo:2007}, and thus omitted. The last statement is obvious.

(b) Let us arbitrarily fix some $x\in S,$ and put
$\hat{\gamma}(\cdot)=\delta_x(\cdot).$ It is obvious that $\hat{\gamma}$ satisfies Condition
\ref{MDPfinitenesscons}(a). Suppose now there is another solution
$v^\ast\in \textbf{B}_{w'}(S)$ to the Bellman equation
(\ref{MDPdpeqn}). But then it follows from part (a) of this theorem that
$\inf_{\pi}V_0(\pi)=u^\ast(x)=v^\ast(x).$

(c) We observe that the Bellman function $u^\ast$ is feasible for linear program
(\ref{MDPunconstrainedDLP}). Consider any function $v$ that is also
feasible for linear program (\ref{MDPunconstrainedDLP}). Therefore, by referring to Lemma \ref{MDPDPinequality}(b), we have that under any Markov policy $\pi,$
$v(x)\le V_0(x,\pi).$ Now suppose $\int_S\gamma(dy)v(y)>\int_S\gamma(dy)u^\ast(y).$ Then there exist some
$\hat{x} \in S$ and constant $\delta>0$ such that $u^\ast(\hat{x})<v(\hat{x})-\delta.$  Hence,
$u^\ast(\hat{x})<V_0(\hat{x},\pi)-\delta,$ where $\pi$ is any Markov policy. But this contradicts part (a) of this
theorem. Therefore, any feasible solution $v$ to linear program (\ref{MDPunconstrainedDLP}) satisfies $\int_S\gamma(dy)v(y)\le \int_S\gamma(dy)u^\ast(y),$ as required.

(d) From part (c) of this theorem, we know that the optimal
value of linear program (\ref{MDPunconstrainedDLP}) is given by $\int_S u^\ast(y)\gamma(dy).$ Therefore, if some feasible solution $v$ to linear program (\ref{MDPunconstrainedDLP}) satisfies
$u^\ast(x)=v(x)$ a.s. with respect to $\gamma$, then it solves the linear program, too. Hence we
conclude the sufficiency part of the statement.

As for the necessity, let $v$ be any optimal solution to linear program (\ref{MDPunconstrainedDLP}). Suppose the relation of $v=u^\ast$ a.s. with respect to $\gamma$ is false. Then there exist measurable subsets $\Gamma_1,\Gamma_2\subseteq S$, such that the following conditions are satisfied: $\Gamma_1\bigcap
\Gamma_2=\emptyset,$  $v(x)>u^\ast(x)$ on $\Gamma_1,$ $v(x)<u^\ast(x)$ on
$\Gamma_2,$ $v(x)=u^*(x)$ on $S\setminus\Gamma_1\setminus\Gamma_2,$ and the  case $\gamma(\Gamma_1)=\gamma(\Gamma_2)=0$ is excluded.
Now let us define a function $\hat{v}$ by
$\hat{v}(x)=I\{x\in S\setminus\Gamma_2\}v(x)+I\{x\in
\Gamma_2\}u^\ast(x),$ which is feasible for linear program (\ref{MDPunconstrainedDLP}).
Indeed, firstly, it is evident that $\hat{v}\in
\textbf{B}_{w'}(S)$. Secondly, we have that $\forall~x\in S\setminus\Gamma_2,$
\begin{eqnarray*}
&&\frac{1}{\alpha}c_0(x,a)-\hat{v}(x)+\frac{1}{\alpha}\int_S\hat{v}(y)q(dy|x,a)\\
&=&\frac{1}{\alpha}c_0(x,a)-v(x)+\frac{1}{\alpha}\int_{S\setminus\Gamma_2}v(y)q(dy|x,a)+\frac{1}{\alpha}\int_{\Gamma_2}u^\ast(y)q(dy|x,a)\\
&\ge&
\frac{1}{\alpha}c_0(x,a)-v(x)+\frac{1}{\alpha}\int_{S\setminus\Gamma_2}v(y)q(dy|x,a)+\frac{1}{\alpha}\int_{\Gamma_2}v(y)q(dy|x,a)\ge
0,
\end{eqnarray*}
and $\forall~x\in \Gamma_2,$
\begin{eqnarray*}
&&\frac{1}{\alpha}c_0(x,a)-\hat{v}(x)+\frac{1}{\alpha}\int_S\hat{v}(y)q(dy|x,a)\\
&=&\frac{1}{\alpha}c_0(x,a)-u^\ast(x)+\frac{1}{\alpha}\int_{S\setminus\Gamma_2}v(y)q(dy|x,a)+\frac{1}{\alpha}\int_{\Gamma_2}u^\ast(y)q(dy|x,a)\\
&\ge&
\frac{1}{\alpha}c_0(x,a)-u^\ast(x)+\frac{1}{\alpha}\int_{S\setminus\Gamma_2}u^\ast(y)q(dy|x,a)+\frac{1}{\alpha}\int_{\Gamma_2}u^\ast(y)q(dy|x,a)\ge
0.
\end{eqnarray*}
However,
$\int_S\hat{v}(y)\gamma(dy)=\int_{S\setminus{\Gamma_2}}v(x)\gamma(dx)+\int_{S\setminus
{\Gamma_2}}u^\ast(x)\gamma(dx)> \int_Sv(x)\gamma(dx),$ which is a
contradiction against that $v$ is optimal for linear program (\ref{MDPunconstrainedDLP}). Now the necessity part follows.
\hfill $\Box$
\bigskip

\noindent\textbf{Proof of Theorem \ref{t6}:}
(a) We take functions $w$ and $w'$ in the form
  $$w(x)=\left\{\begin{array}{ll}
  1, & \mbox{ if } x=0; \\
  \frac{1}{x^4}, & \mbox{ if } x\in(0,1];
  \end{array}\right. $$
$$w'(x)=\left\{\begin{array}{ll}
  1, & \mbox{ if } x=0; \\
  \frac{1}{x^2}, & \mbox{ if } x\in(0,1],
  \end{array}\right. $$
and put $S_0=\{0\}$, $S_l=S_0\cup \left(\frac{1}{l+1},1\right]$, $l=1,2,\dots.$ Now Condition \ref{MDPRegularitycons}(a,c) is obviously satisfied.

Condition \ref{MDPRegularitycons}(b) can be verified for $\rho\defi 4\lambda$ and $b=0$ as follows:\\
-- if $x=0$ then
  $$\int_S q(dy|x,a) w(y)=5\lambda\int_0^1\frac{1}{y^4} y^4 dy-\lambda=4\lambda=\rho w(0);$$
-- if $x\in(0,1]$ then
  $$\int_S q(dy|x,a)w(y)=\frac{a}{x} w(0)-\frac{a}{x} w(x)=\frac{a}{x}\left(1-\frac{1}{x^4}\right)\le 0<\rho w(x).$$

For Condition \ref{MDPfinitenesscons}, it is sufficient to notice that $\forall~x\in (0,1],$
$$\inf_{a\in A(x)} c_0(x,a)=\left\{\begin{array}{ll}
C_1 x-\frac{1}{4C_2 x^2}, & \mbox{ if } \frac{1}{2C_2}<\bar A; \\ \ \\
C_1 x+C_2\frac{\bar A^2}{x^2}-\frac{\bar A}{x^2}, & \mbox{ otherwise},
\end{array}\right. $$
$\inf_{a\in A(0)}c_0(0,a)=0,$ and  $\alpha>4\lambda=\rho.$

Condition \ref{MDPmeas1} and Condition \ref{MDPstrongregularitycons} are trivially satisfied because
$$q_x(a)=\left\{\begin{array}{ll}
\lambda, & \mbox{ if } x=0, \\ \ \\
\frac{a}{x}, & \mbox{ if } x\in(0,1],
  \end{array}\right. $$
$\forall~x\in(0,1],A(x)=[0,\frac{\bar A}{x}]$, and $A(0)=\{0\}.$

Condition \ref{MDPdpcons}(b,c,d) can be verified similarly to what is presented above by taking $\rho'=\frac{2\lambda}{3}$, $b'=0$. Since $\forall~x\in(0,1],\bar q_x\le \frac{\bar A}{x^2}$ and $\bar{q}_{0}=\lambda,$ Condition \ref{MDPdpcons}(a) is also satisfied.

Finally, Condition \ref{MDPmeas2} obviously holds.

(b) If we denote $z^{(n+1)}=f(z^{(n)})$ then, for $z>\frac{\epsilon}{2}>0$, where $\epsilon>0$ is any fixed constant, function $f$ is differentiable:
 $$\frac{df}{dz}=\frac{-5\lambda}{\alpha+\lambda}\int_0^1\frac{\partial u(y,z)}{\partial z} y^4 dy,$$
 where
\begin{eqnarray*}
\frac{\partial u(x,z)}{\partial z}&=&-1+\frac{\alpha C_2 x^2}{\sqrt{\alpha^2 C_2^2 x^4+C_1C_2 x^3+\alpha C_2 x^2z}}\\
&=&\frac{\alpha C_2 x^2-\sqrt{\alpha^2 C_2^2 x^4+C_1C_2 x^3+\alpha C_2 x^2z}}{\sqrt{\alpha^2 C_2^2 x^4+C_1C_2 x^3+\alpha C_2 x^2z}}\in (-1,0),\forall~x\in(0,1],
\end{eqnarray*}
so that $\forall~z\in(\frac{\epsilon}{2},\infty), 0<\frac{df}{dz}<\frac{\lambda}{\alpha+\lambda}<1$.

It remains to estimate $z^{(1)}$:
\begin{eqnarray*}
u^{(1)}(x) = -2\alpha C_2 x^2+2\sqrt{\alpha^2 C_2^2 x^4+C_1C_2 x^3}\le -2\alpha C_2 x^2+\left(2\alpha C_2 x^2+\frac{C_1 x}{\alpha}\right)=\frac{C_1 x}{\alpha},\forall~x\in(0,1];
\end{eqnarray*}
\begin{eqnarray*}
z^{(1)} \ge 1-\frac{5\lambda C_1}{\alpha(\alpha+\lambda)}\int_0^1 y dy>1-\frac{C_1}{2\alpha}\ge 0
\end{eqnarray*}
because $\alpha>4\lambda$ and $C_1<2\alpha$. The map $z\to f(z)$ is contracting on $[\epsilon,\infty)$, e.g., for $\epsilon=z^{(1)}$.
Since
\begin{eqnarray*}
f\left(\frac{10}{7} C_2\lambda+\frac{\alpha+\lambda}{\alpha}\right) < 1+\frac{5\lambda}{\alpha+\lambda}\left[\int_0^1 \left(2\alpha C_2 x^2+\frac{10}{7} C_2\lambda +\frac{\alpha+\lambda}{\alpha}\right) x^4 dx\right]=\frac{10}{7} C_2\lambda+\frac{\alpha+\lambda}{\alpha},
\end{eqnarray*}
we conclude that $z^*< \frac{10}{7} C_2\lambda+\frac{\alpha+\lambda}{\alpha}$.

(c) Clearly, function $u^*(x)$ (supplemented by $u^\ast(0)=1-z^\ast$) is bounded; hence $u^*\in \textbf{B}_{w'}(S)$. Therefore, according to Theorem \ref{MDPbellmanthm}, it is sufficient to check that $u^*$ solves equation (\ref{MDPdpeqn}) and $\phi^*$ provides the infimum.

Expression in the parenthesis of (\ref{MDPdpeqn}) equals
$$\lambda\int_0^1 u^*(y) 5 y^4 dy-\lambda u^*(0) \mbox{ if } x=0, $$
and
$$C_1x+C_2 a^2-\frac{a}{x} +\frac{a}{x} u^*(0)-\frac{a}{x} u^*(x) \mbox{ if } x\in(0,1].$$
Therefore,
  $$u^*(0)=\frac{5\lambda}{\alpha+\lambda}\int_0^1 u^*(y) y^4 dy$$
  and $\phi^*(x)$ given by (\ref{MDPexamplepolicy}) provides the infimum. (Note that $u^*(x)+z^*\ge -2\alpha C_2 x^2+2\sqrt{\alpha^2 C_2^2 x^4}=0.$)
Finally, at $x>0$, the RHS of (\ref{MDPdpeqn}) equals $C_1 x-\frac{(u^*(x)+z^*)^2}{4x^2 C_2}$, and equation
  $$4\alpha C_2 x^2 u^*(x)=4C_1C_2 x^3-(u^*(x))^2-2 u^*(x) z^*-(z^*)^2$$
holds because
  $$u^*(x)=-2\alpha C_2 x^2-z^*+2\sqrt{\alpha^2 C_2^2 x^4+C_1C_2 x^3+\alpha C_2 x^2 z^*}.$$
\hfill $\Box$
\vfill


\end{document}